\documentclass[twoside,11pt]{amsart}
\usepackage{graphicx, lscape}
\usepackage{amsmath, amsfonts, amssymb, ifthen}

\newcommand{\hide}[1]{} 

\renewcommand{\(}{\left(}
\renewcommand{\)}{\right)}  
\newcommand{\abs}[1]{\left\arrowvert #1\right\arrowvert}

\newcommand{\field}[1]{\mathbb{#1}}
\newcommand{\CC}{\ensuremath{\field{C}}} 
\newcommand{\RR}{\ensuremath{\field{R}}} 
\newcommand{\NN}{\ensuremath{\field{N}}} 
\newcommand{\QQ}{\ensuremath{\field{Q}}}
\newcommand{\DD}{\ensuremath{\field{D}}} 
\newcommand{\ZZ}{\ensuremath{\field{Z}}}

\newcommand{\An}[2]{\ensuremath{\field{A}}\left(#1, #2\right)} 

\newcommand{\Arg}[1]{\text{Arg}(#1)}
\newcommand{\HH}{\mathcal{H}}

\newcommand{\Ln}{\text{Ln}}

 \newcommand{\cDD}{\overline{\DD}}
 
 \newcommand{\Pnc}{P_{n,c}}
 \newcommand{\MPnc}{M_n(P)}
 
 \newcommand{\JPnc}{J(P_{n,c})}
 \newcommand{\KPnc}{K(P_{n,c})}
 
   \newcommand{\Ftc}{F_{t,c}}
 \newcommand{\MFtc}{M_t(F)}
 
 \newcommand{\JFtc}{J(F_{t,c})}
 \newcommand{\KFtc}{K(F_{t,c})}
 
\newcommand{\Rna}{R_{n, a}}
 \newcommand{\MRna}{M_n(R_0)}
 
 \newcommand{\JRna}{J(R_{n,a})}
 \newcommand{\KRna}{K(R_{n,a})}

\newcommand{\Rnca}{R_{n,c,a}}
 \newcommand{\MRnca}{M_n(R_c)}
 
 \newcommand{\JRnca}{J(R_{n,c,a})}
 \newcommand{\KRnca}{K(R_{n,c,a})}
\newcommand{\MtwoRnca}{M^2_n(R_c)}
\newcommand{\MoneRnca}{M^1_n(R_c)}

\newcommand{\eps}{\epsilon}

\theoremstyle{plain}
\newtheorem{thm}{Theorem}[section]
\newtheorem{cor}[thm]{Corollary}
\newtheorem{lem}[thm]{Lemma}
\newtheorem{prop}[thm]{Proposition}

\theoremstyle{definition}
\newtheorem{defn}[thm]{Definition}

\newcommand{\showcomments}{yes}

\newsavebox{\commentbox}
\newenvironment{comment}%
{\ifthenelse{\equal{\showcomments}{yes}}%
{\footnotemark
    \begin{lrbox}{\commentbox}
    \begin{minipage}[t]{1.25in}\raggedright\sffamily\small
    \footnotemark[\arabic{footnote}]}
{\begin{lrbox}{\commentbox}}}%
{\ifthenelse{\equal{\showcomments}{yes}}%
{\end{minipage}\end{lrbox}\marginpar{\usebox{\commentbox}}}
{\end{lrbox}}}

\newcommand{\drawfigcauliflower}{\includegraphics[height=6cm]{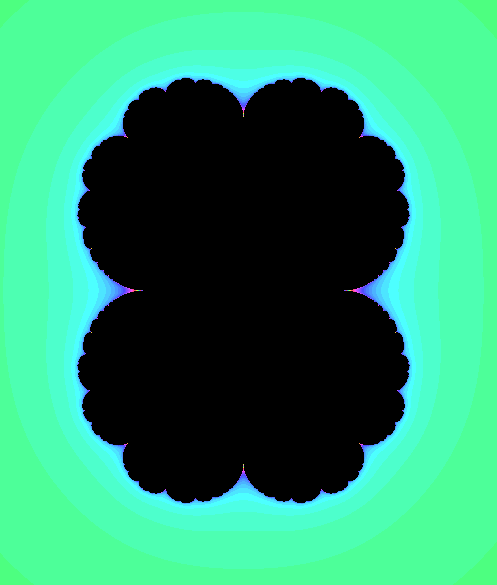}}
\newcommand{\drawfigcaulicantor}{\includegraphics[height=6cm]{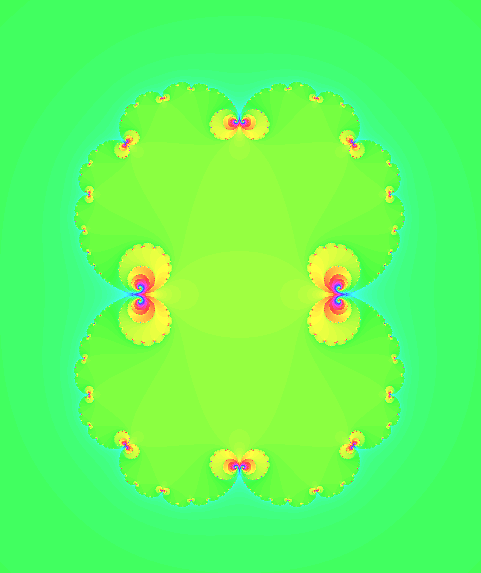}}

\newcommand{\drawfigMn }{{\includegraphics[height=4cm]{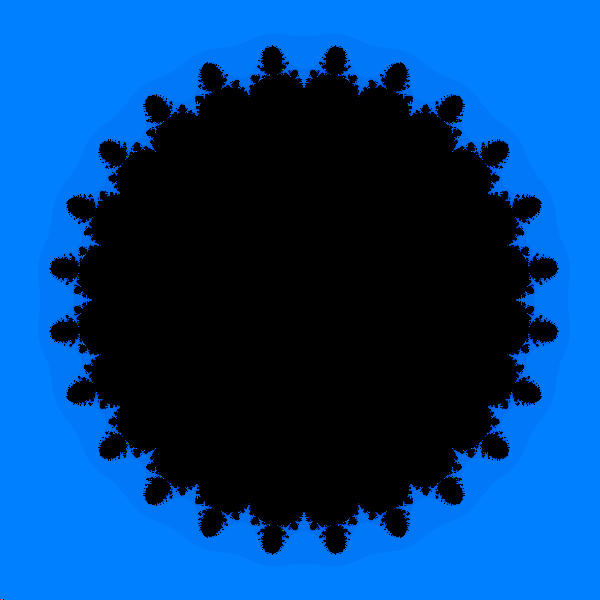}} }
\newcommand{\drawfigJDn }{ {\includegraphics[height=4cm]{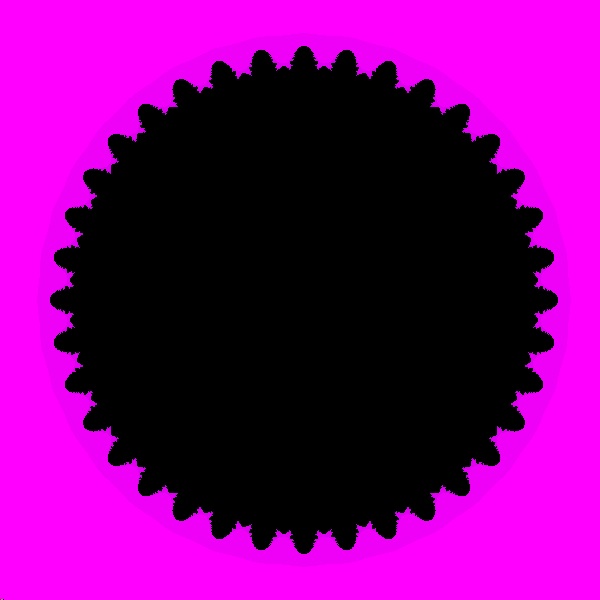}} }
\newcommand{\drawfigJSn }{ {\includegraphics[height=4cm]{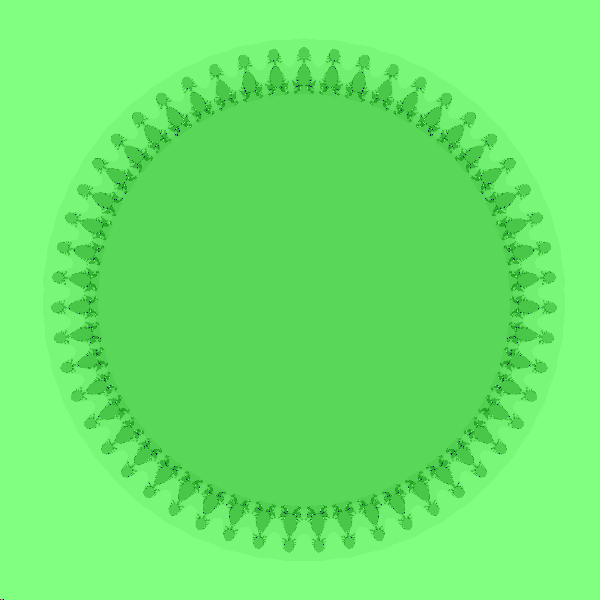}} }

\newcommand{ \drawfigJRnca }{{\includegraphics[height=6cm]{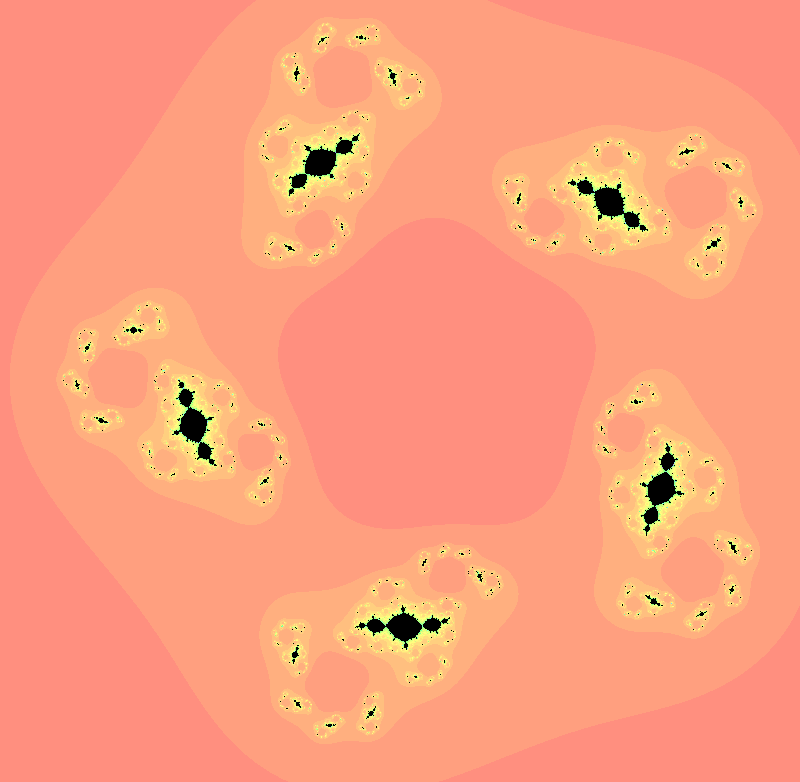}}  }
\newcommand{ \drawfigMRnca }{ {\includegraphics[height=6cm]{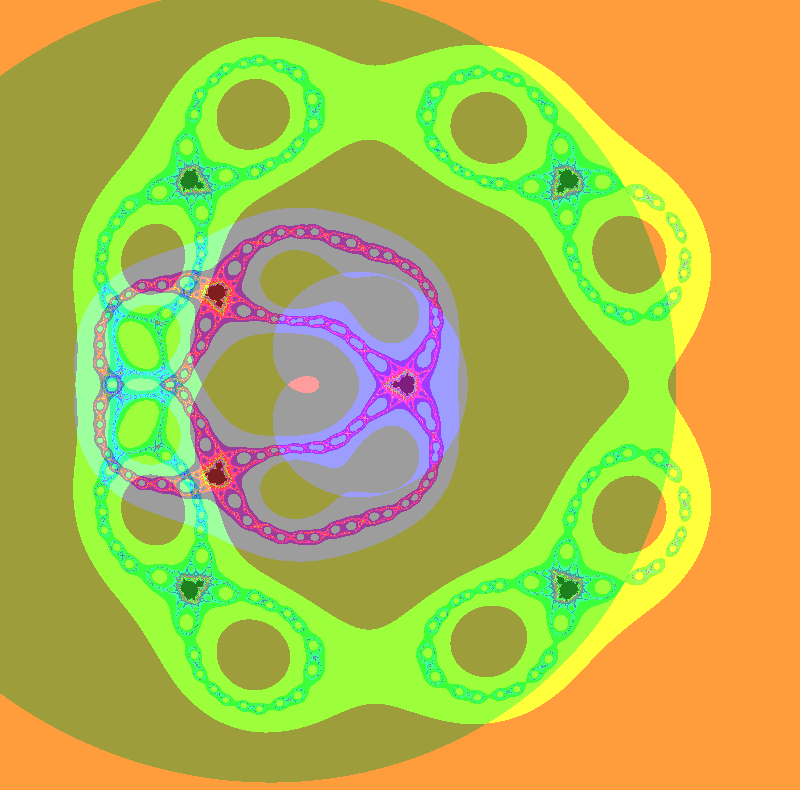}} }

\newcommand{ \drawfigKPoutnlow }{ {\includegraphics[height=5cm]{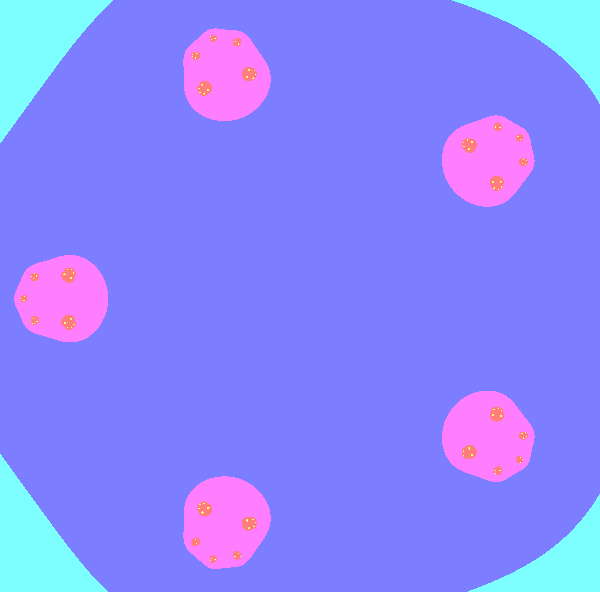}}}
\newcommand{ \drawfigKPoutnhi }{ {\includegraphics[height=5cm]{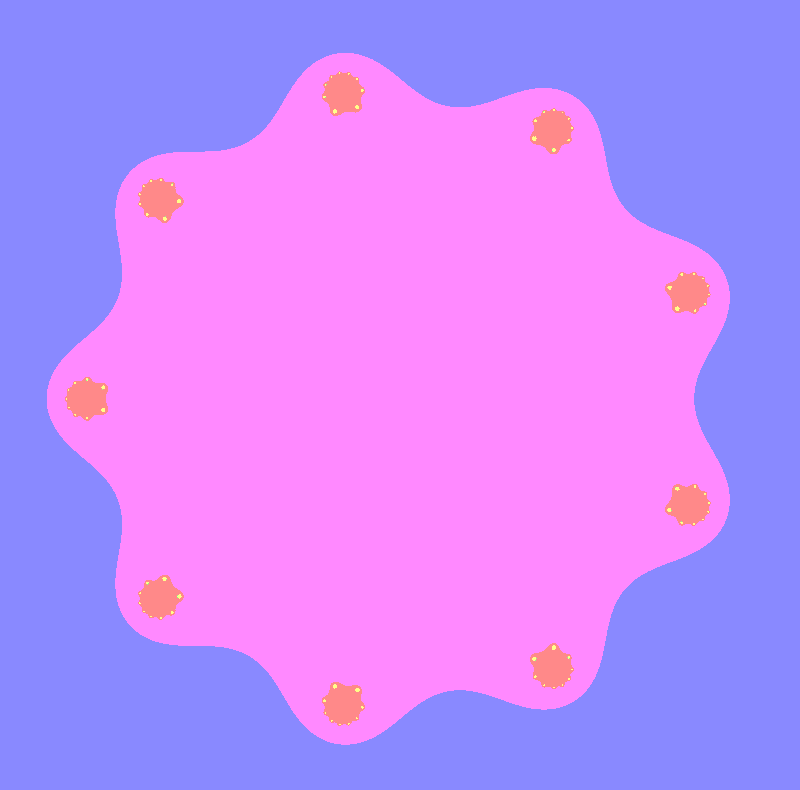}}  }

\newcommand{ \drawfigKPinnlow }{ {\includegraphics[height=5cm]{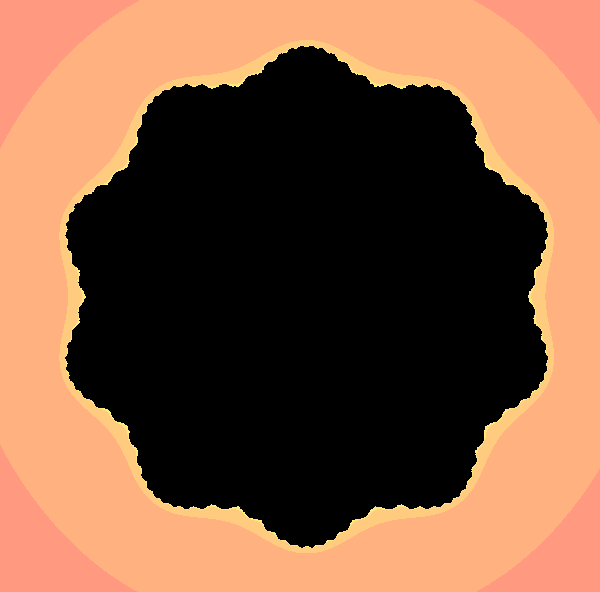}}  }
\newcommand{ \drawfigKPinnhi }{ {\includegraphics[height=5cm]{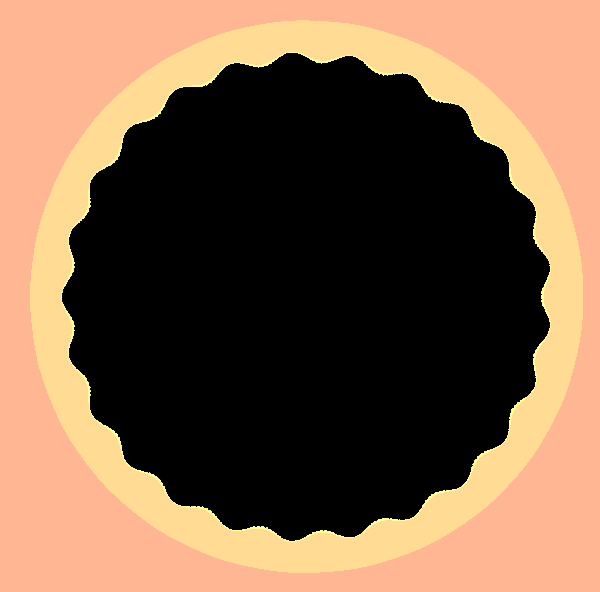}} }

\newcommand{ \drawfigMPnlo }{ {\includegraphics[height=4cm]{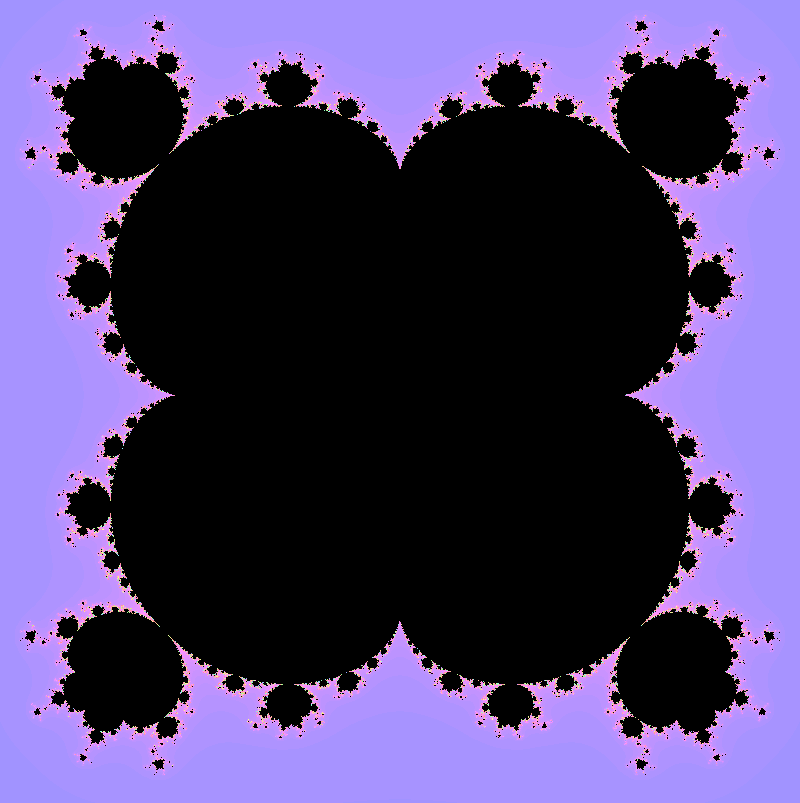}}} 
\newcommand{ \drawfigMPnmed }{ {\includegraphics[height=4cm]{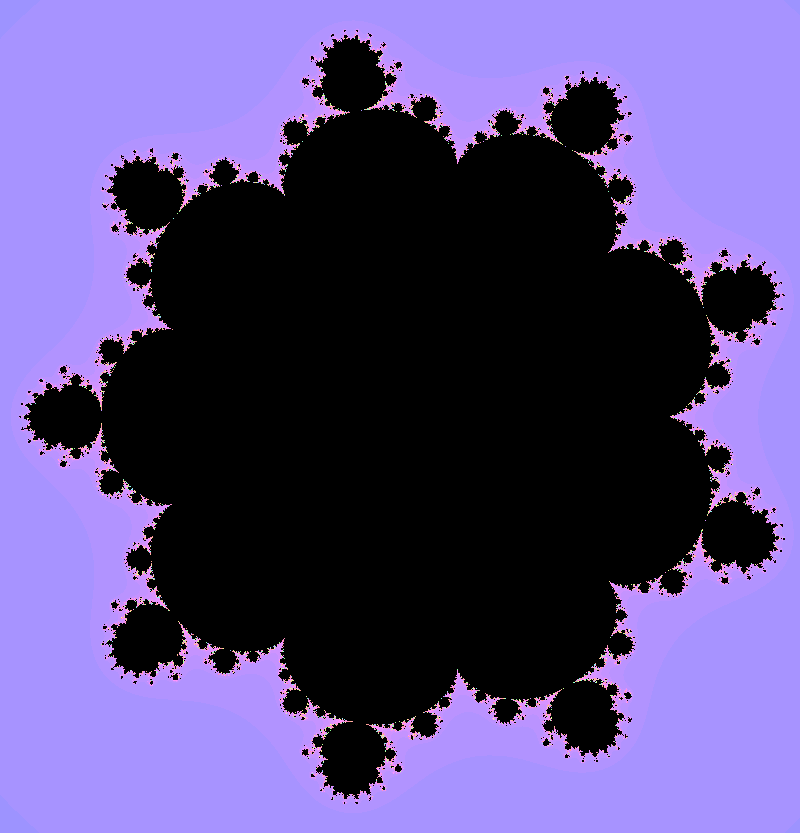}}  }
\newcommand{ \drawfigMPnhi }{ {\includegraphics[height=4cm]{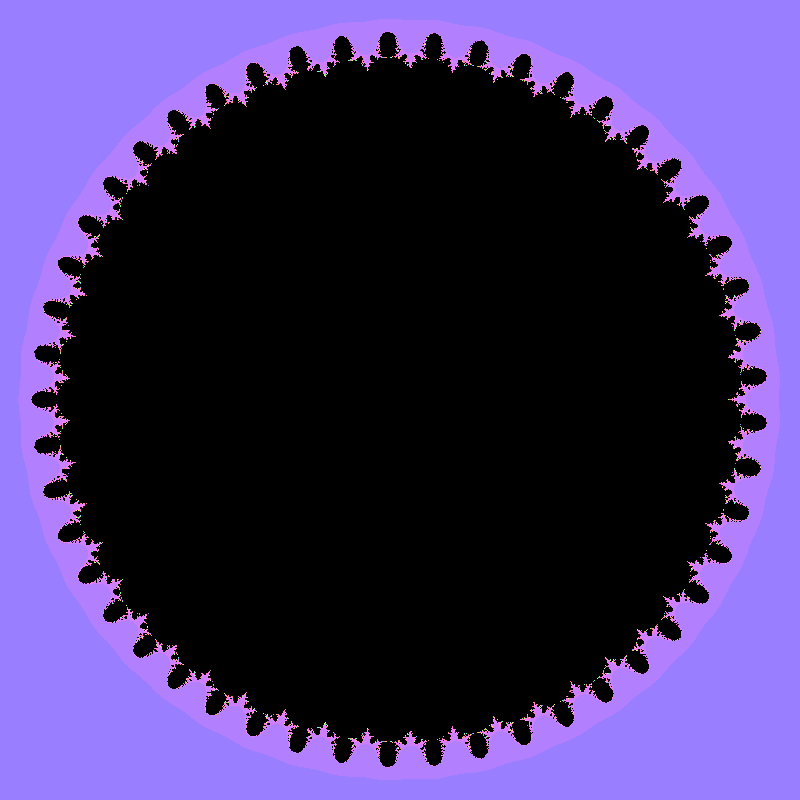}}  }

\newcommand{ \drawfigKPincone }{ {\includegraphics[height=6cm]{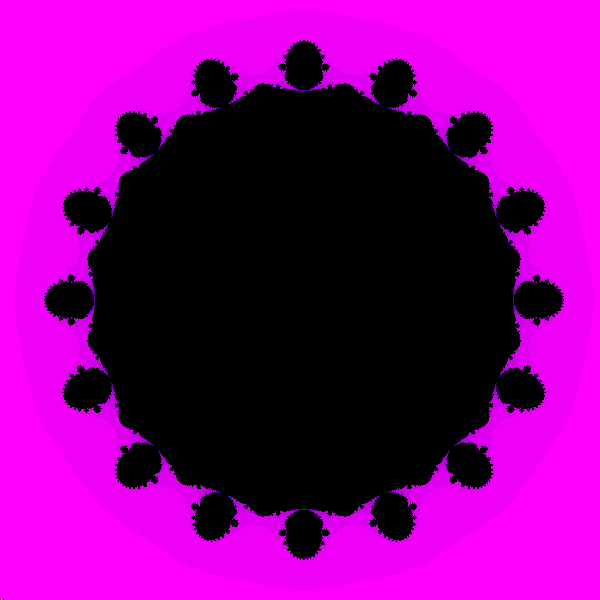}}  }
\newcommand{ \drawfigKPoutcone }{  {\includegraphics[height=6cm]{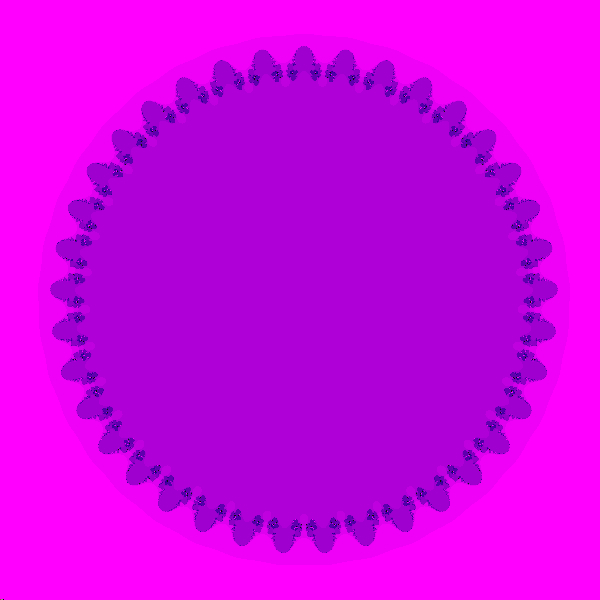}}  }

\newcommand{ \drawfigKFtcin }{ {\includegraphics[height=5.5cm]{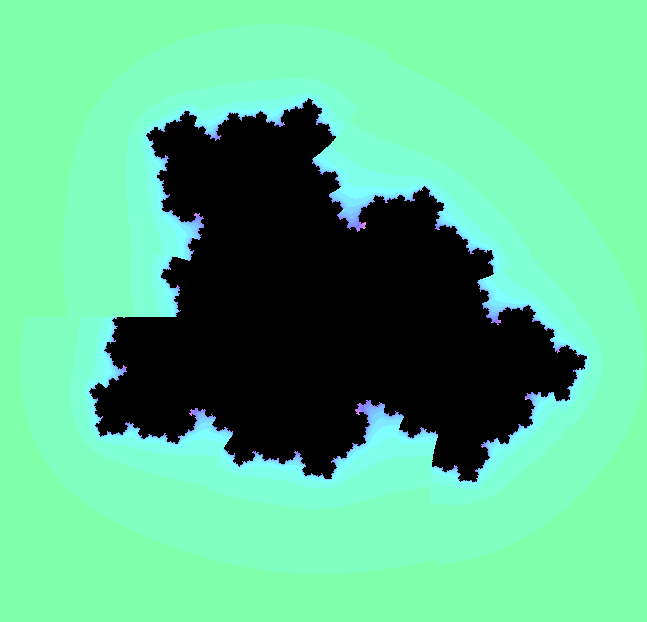}}  }
\newcommand{\drawfigKFtcout }{ {\includegraphics[height=5.5cm]{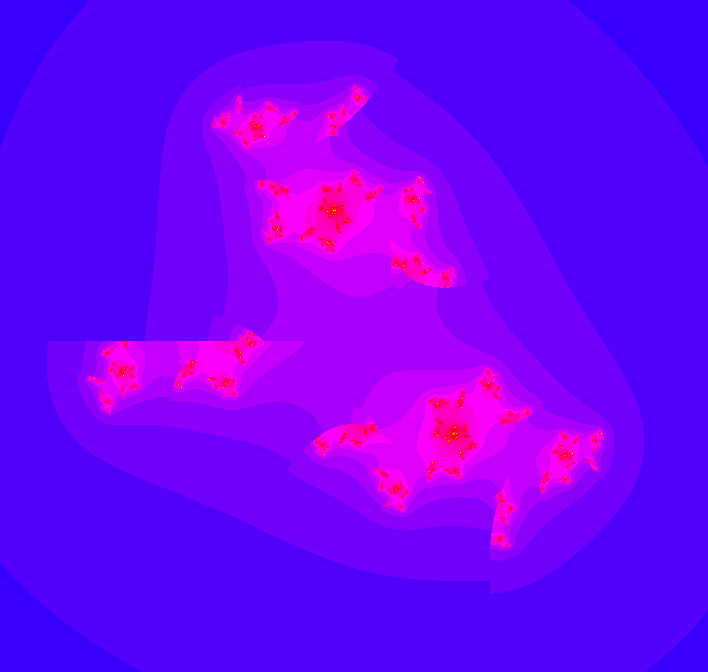}} }

\newcommand{\drawfigMFtcnlo}{ {\includegraphics[height=6cm]{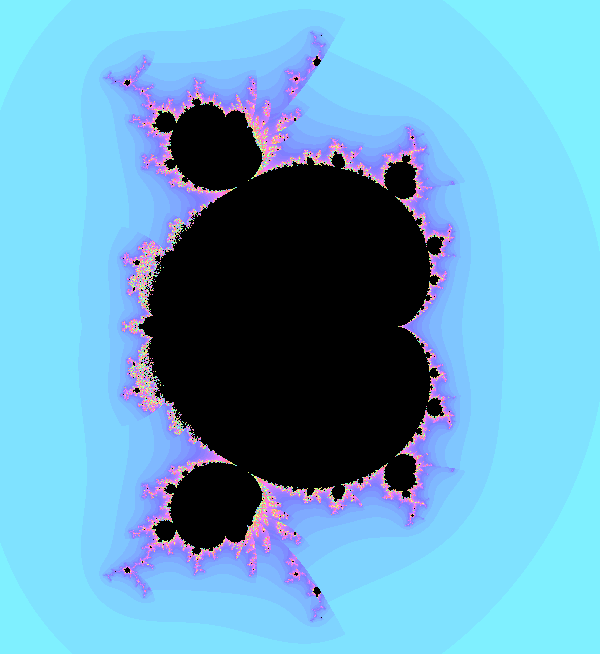}} }
\newcommand{ \drawfigMFtcnhi }{ {\includegraphics[height=6cm]{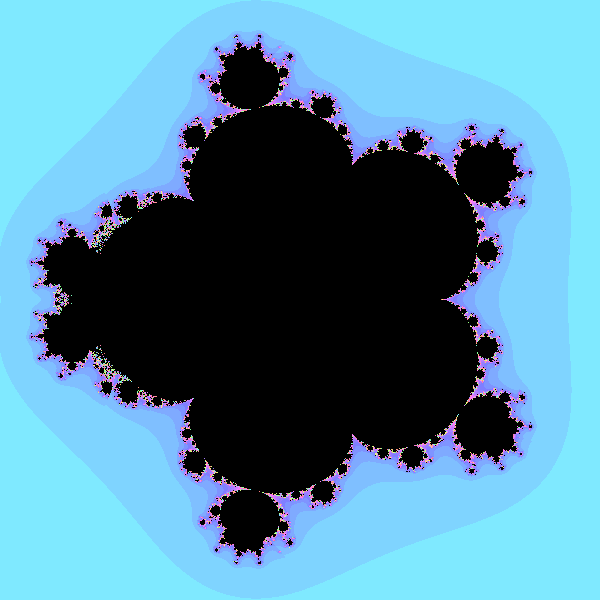}}  }

\newcommand{ \drawfigKRnaexmp }{ {\includegraphics[height=6cm]{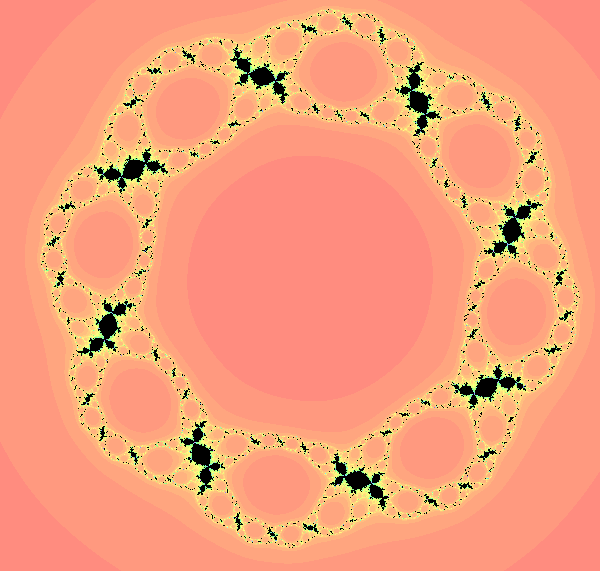}} }
\newcommand{ \drawfigKRncaexmp }{ {\includegraphics[height=6cm]{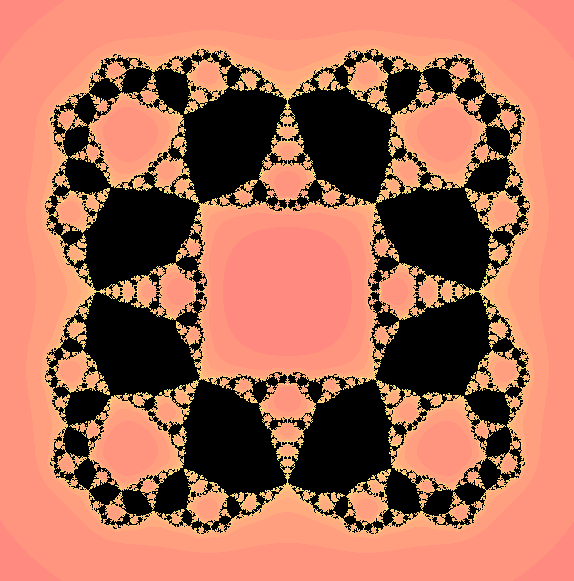}} }

\newcommand{ \drawfigMRnaexmpnlo }{ {\includegraphics[height=6cm]{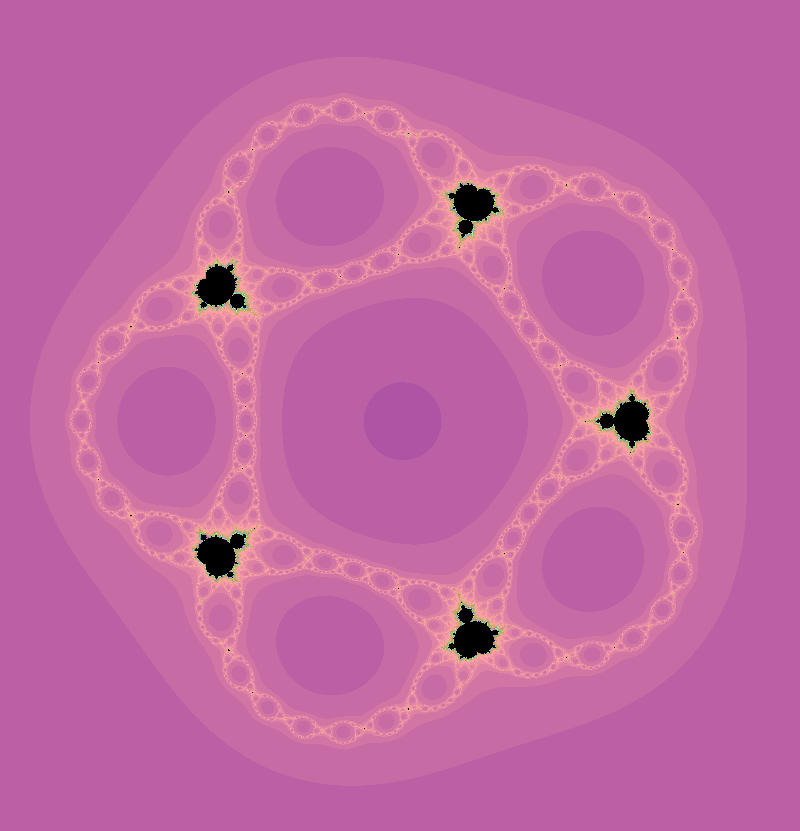}}  }
\newcommand{ \drawfigMRncaexmpnhi }{ {\includegraphics[height=6cm]{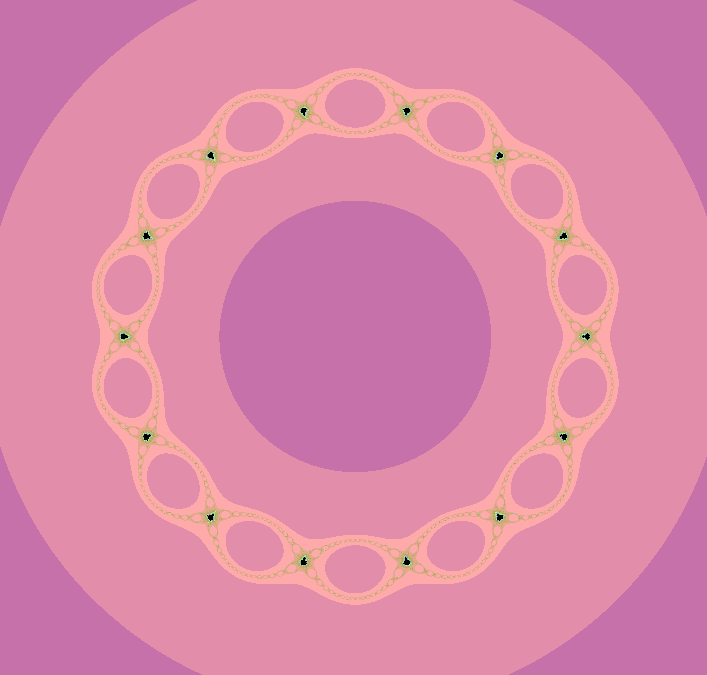}}  }

\newcommand{ \drawfiglimclo }{  {\includegraphics[height=4cm]{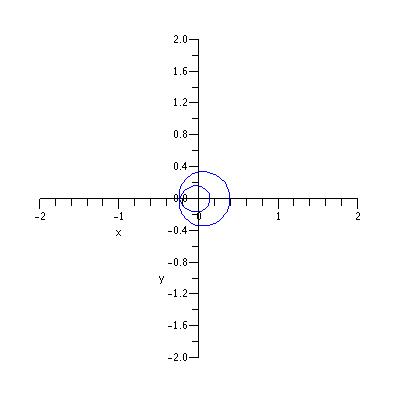}} }
\newcommand{ \drawfiglimcmed }{ {\includegraphics[height=4cm]{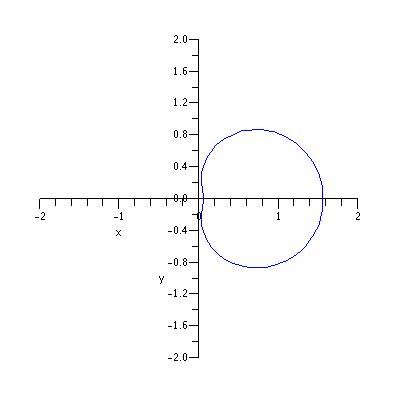}}  }
\newcommand{ \drawfiglimchi }{ {\includegraphics[height=4cm]{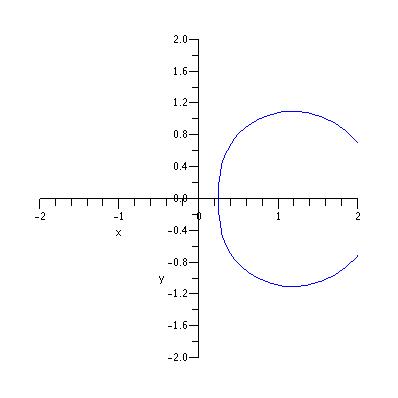}} }

\newcommand{\drawfigMlimclo }{  {\includegraphics[height=4cm]{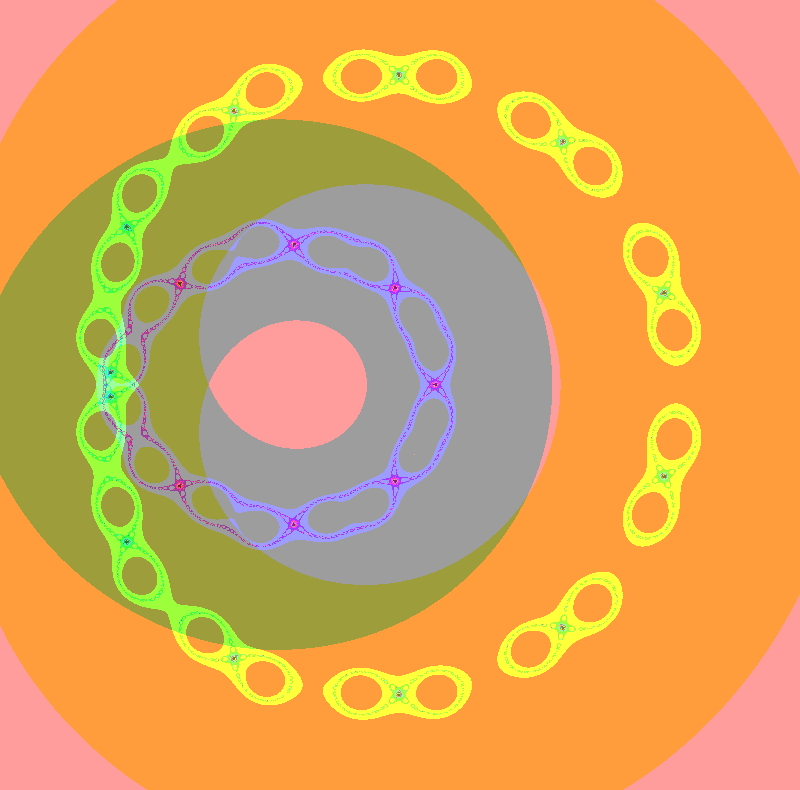}} }
\newcommand{ \drawfigMlimcmed }{  {\includegraphics[height=4cm]{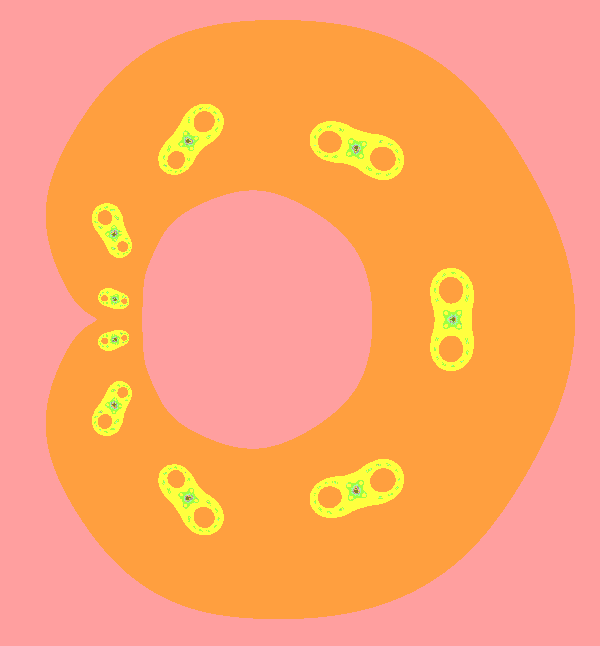}} }
\newcommand{\drawfigMlimchi }{ {\includegraphics[height=4cm]{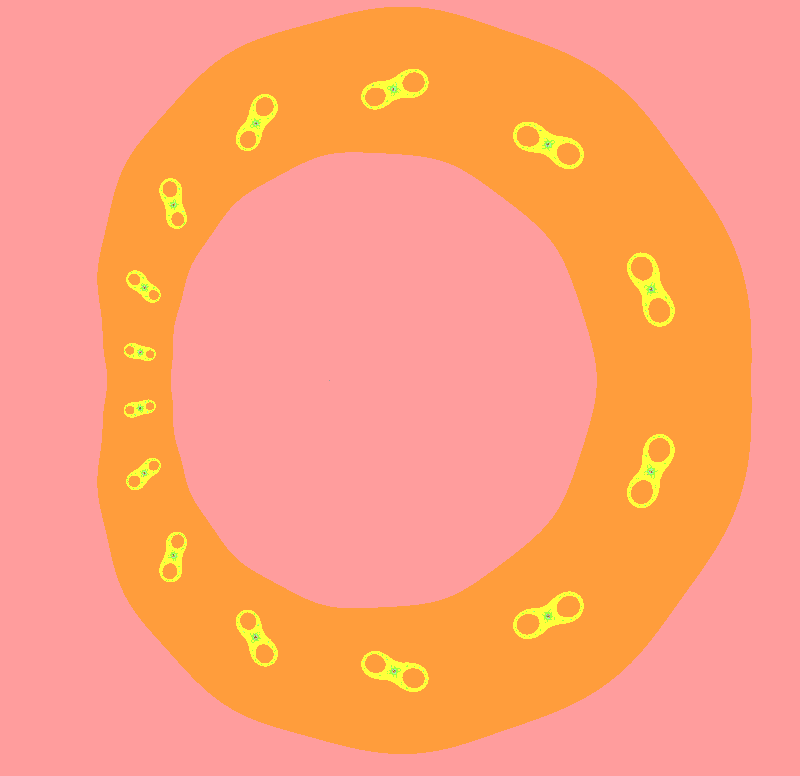}} }

\newcommand{ \drawfigcardrot }{ \includegraphics[height=12cm]{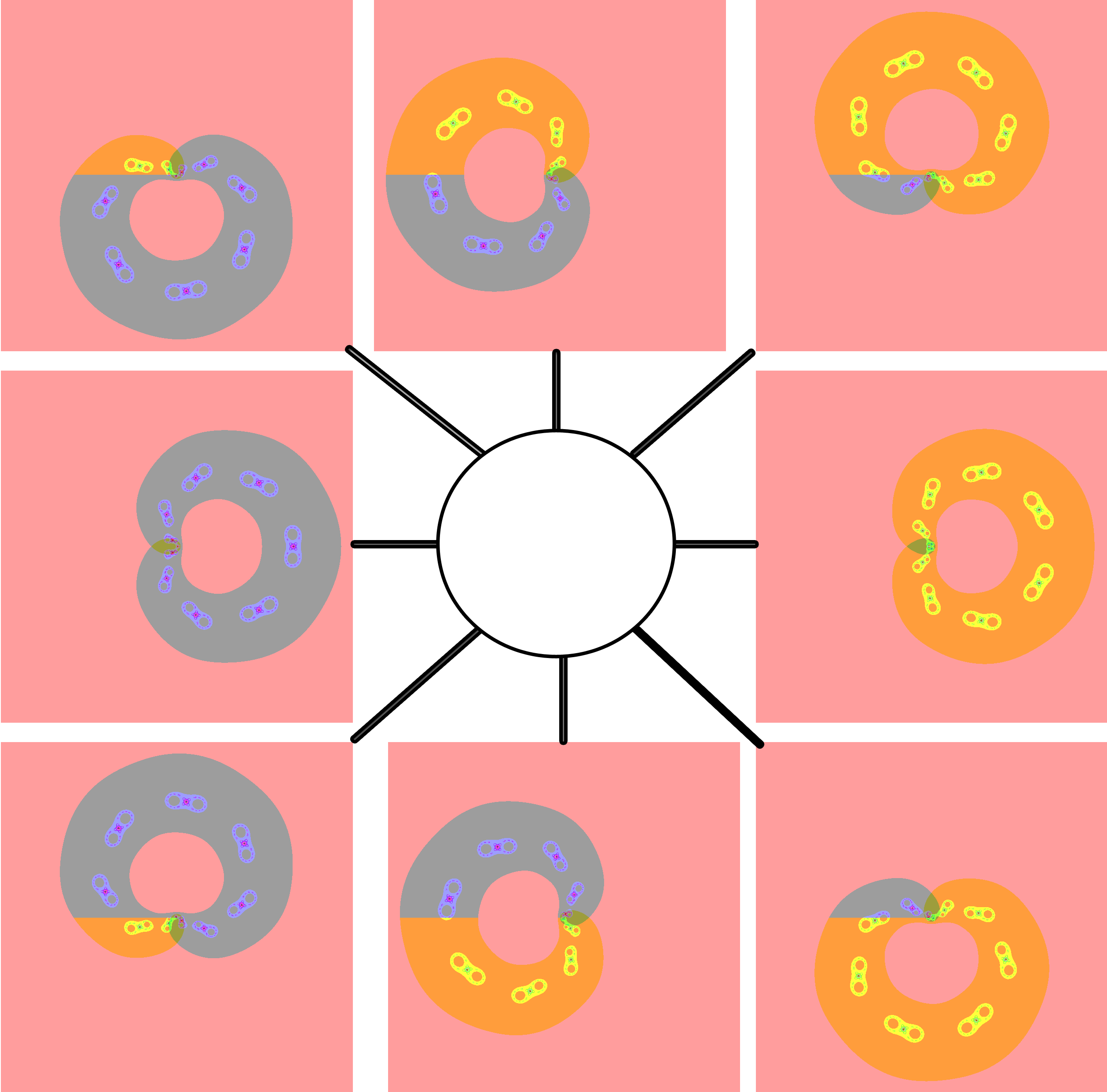}}

\newcommand{ \drawfigfghone }{ \includegraphics[height=7.5cm]{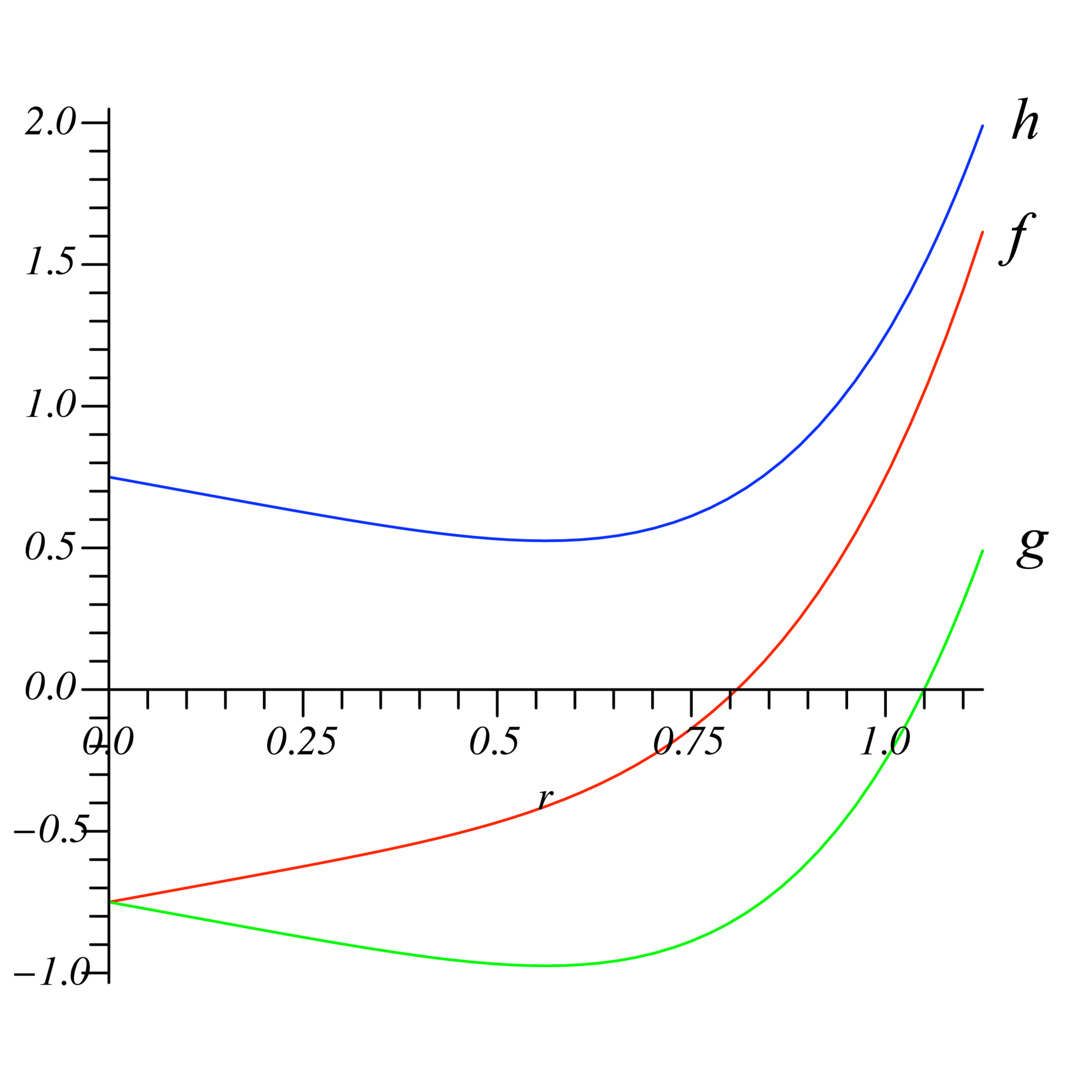}}

\newcommand{\drawfigfghtwo}{ \includegraphics[height=7.5cm]{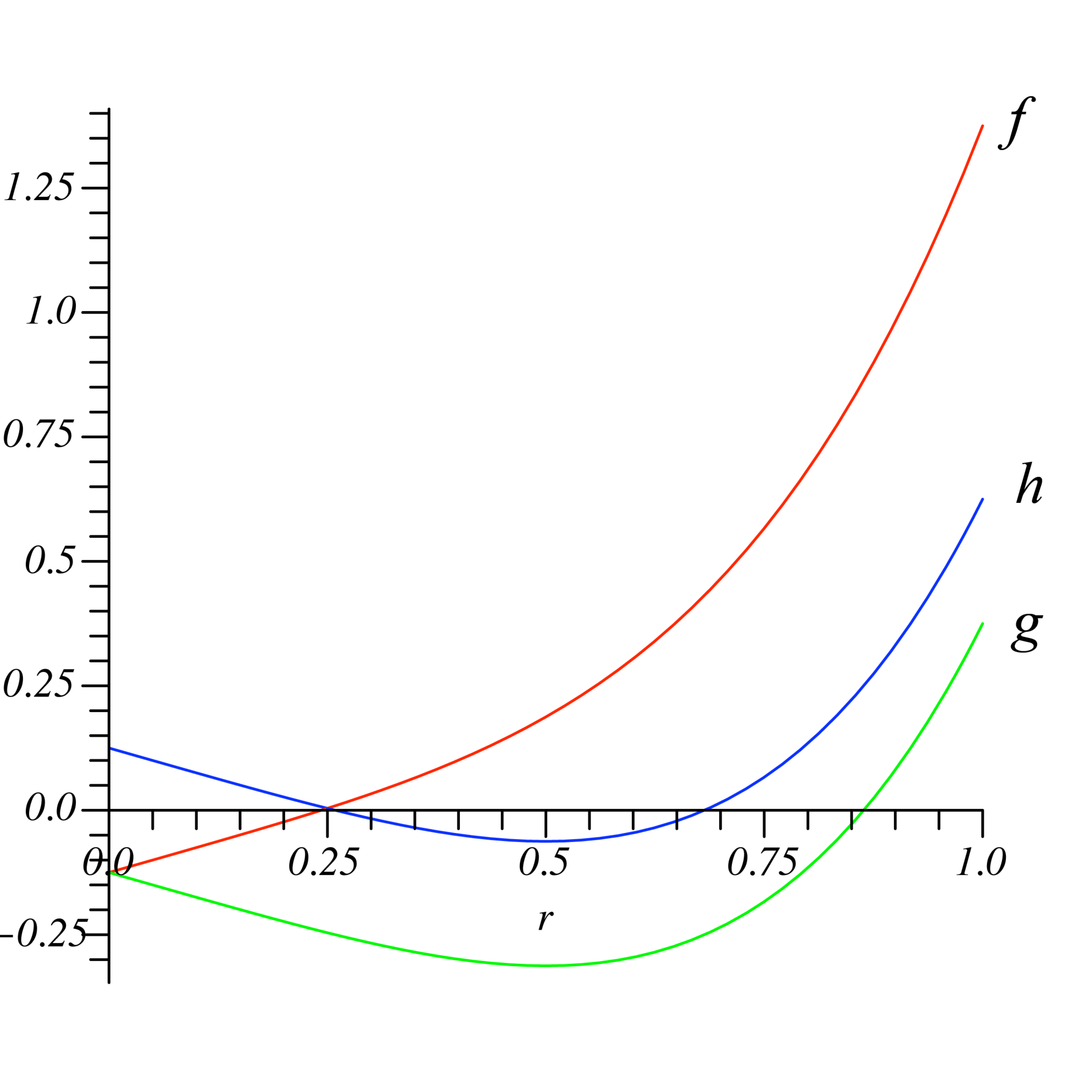}}
\newcommand{\drawfigfghthree }{ \includegraphics[height=5cm]{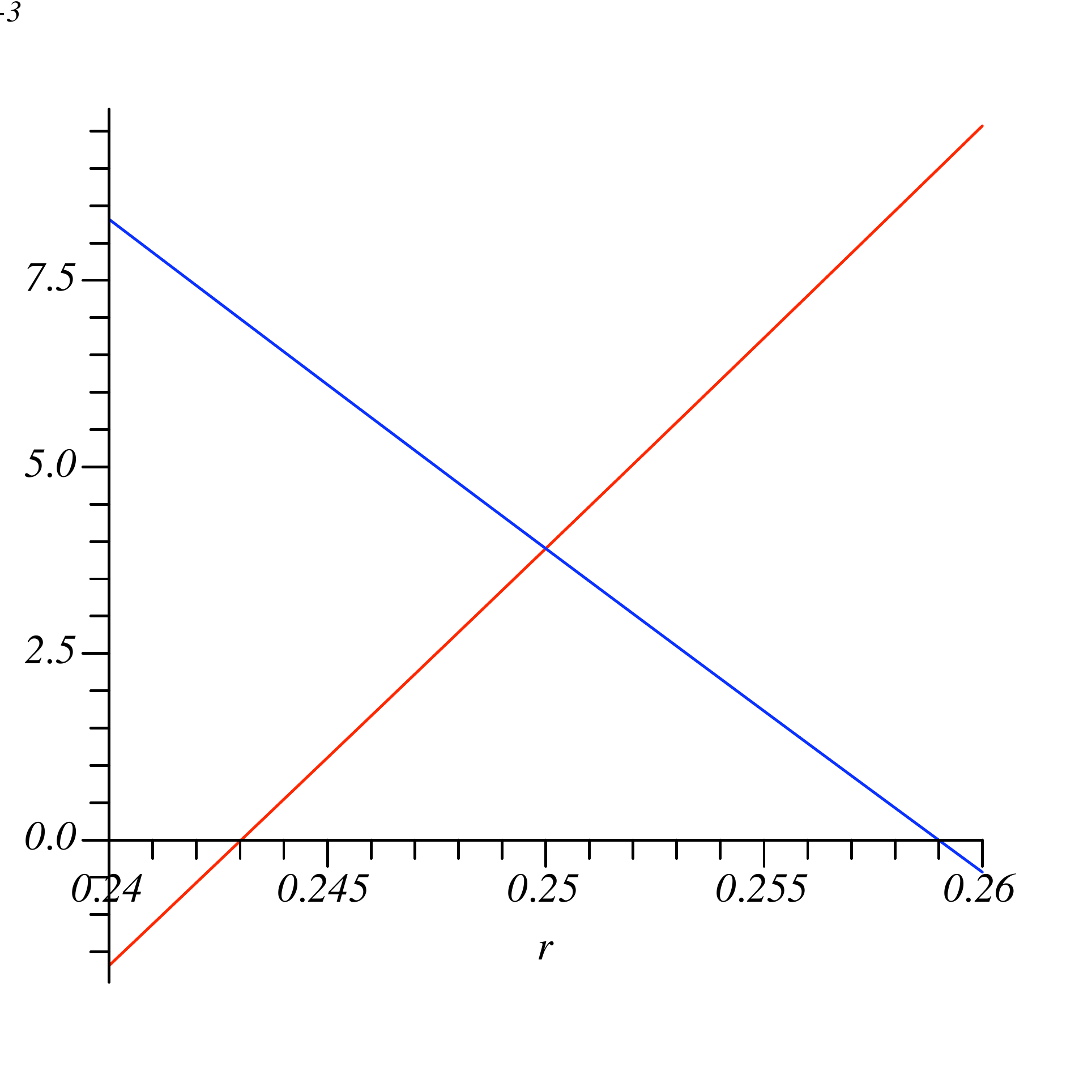}}

\begin{document}

\title{Geometric Limits of Mandelbrot and Julia Sets under degree growth}

\author[S.H. ~Boyd]{Suzanne Hruska Boyd}
\address{Department of Mathematical Sciences\\
University of Wisconsin Milwaukee\\
PO Box 413\\
Milwaukee, WI 53201\\
USA}
\email{sboyd@uwm.edu}
\author[M.J. ~Schulz]{Michael J. Schulz}
\address{Department of Mathematical Sciences\\
University of Wisconsin Milwaukee\\
PO Box 413\\
Milwaukee, WI 53201\\
USA}
\email{mjschulz@uwm.edu}

\date{\today}

\begin{abstract}
First,  for
the family $\Pnc(z) = z^n + c$, we show that the geometric limit of the Mandelbrot sets $\MPnc$ as $n \to \infty$ exists and is the closed unit disk, and that the geometric limit of 
the Julia sets $\JPnc$ as $n$ tends to infinity is the unit circle, at least when $|c| \neq 1$. Then we establish similar results  for some generalizations of this family; namely, the maps $z \mapsto z^t+c$ for real $t \geq 2$ and the rational maps $z  \mapsto z^n + c + a/z^n$.
\end{abstract}

\maketitle

\markboth{\textsc{S. Boyd and M. Schulz}}
  {\textit{}}



\section{Introduction}
\label{sec:introduction}

We begin by considering the family of polynomial maps 
$$\Pnc (z) = z^n + c,$$ where $n \geq 2$ is an integer and $c$ is complex.
These maps are all unicritical (the only critical point is $0$),
so they provide a ``one-step'' generalization of the well known quadratic family $z \mapsto z^2 + c$. 

As usual in polynomial dynamics, for a given $n,c$ we may define the \textit{filled Julia set}  as $K(\Pnc) = \{ z \colon \{\Pnc^m(z)\}_{m=1}^{\infty} $ is bounded$\}$ and the \textit{Julia set} as the topological boundary, $J(\Pnc) = \partial K(\Pnc)$.  For each $n \geq 2$, we define the
\textit{Mandelbrot set} as $\MPnc= \{ c  \colon \{\Pnc^m(0)\}_{m=1}^{\infty} $ is bounded $\}$. The Mandelbrot sets for $n > 2$ are often called multibrot sets.
One may find some basic results on $\Pnc$ in \cite{Dierk}, providing analogs of many fundamental results from the quadratic case to the general degree case.  

\begin{figure}
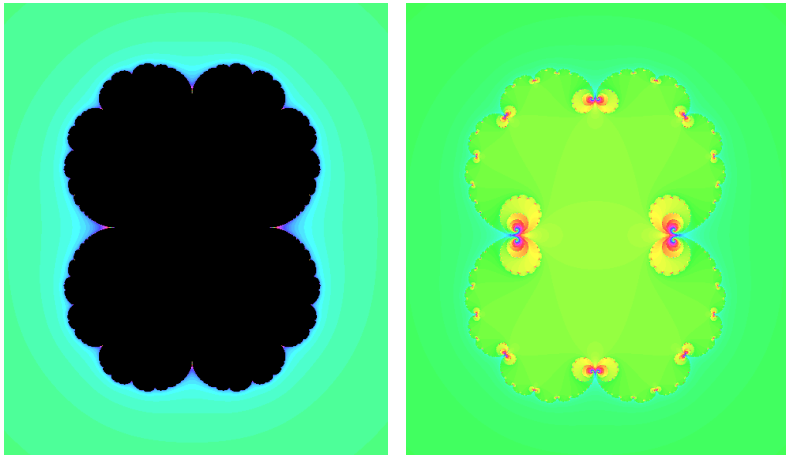

\drawfigcauliflower \ \ 
\drawfigcaulicantor
\caption{\label{fig:cvaries} $\KPnc$ for $n=2$, with $c=1/4$ (L) and $c=.251$ (R). }
\end{figure}

Much study has gone into the exploration of how the dynamics as well as the geometry of the Julia set changes as the parameter $c$ varies.   The geometric limit at $P_{2, 1/4} (z) = z^2 + 1/4$ well illustrates how, although $K_c$ varies upper semi-continuously and $J_c$ varies lower semi-continuously, at $c=1/4$ neither set varies continuously (see Figure~\ref{fig:cvaries}).  This type of ``parabolic implosion'' has been developed by A.\ Douady, A.\ Epstein, J.\ Hubbard, R.\ Oudkerk, P. Lavaurs, M.\ Shishikura, and others, and is still under study.   We also refer the reader to McMullen \cite{Curtbook} for a good introduction to the of stability of polynomial dynamics through small changes in the parameters.  

Here, we take a different approach, and scrutinize these families of dynamical systems $P_{n,c}$ for large $n$.  We look for qualitative similarities among these maps as $n$ tends to infinity.  Of course in the limit, there is no dynamics, no reasonable map $P_{\infty, c}$.  However, it turns out the geometric sets $\KPnc$, $\JPnc$, and even $\MPnc$ have limits that exist as $n\to \infty$.  For $\MPnc$, this limit is simply the closed unit disk.

\begin{thm}  \label{thm:mainMPnc}
Let ${CS}(\hat{\CC})$ denote the collection of all compact subsets of $\hat{\CC}$, then 
\[\lim_{n\to\infty}\MPnc=\cDD.
\]
 under the Hausdorff metric $d_\HH$ in ${CS}(\hat{\CC})$.
\end{thm}

In Section~\ref{sec:background} we give a precise definition of convergence of compact sets.

For comparison purposes, we note that in \cite{DevExp}, it is shown that for the family
 $$Q_{n,\lambda} (z) = \lambda \( 1 + \frac{z}{n} \)^n,$$
 the limit of the Mandelbrot sets $M_n(Q)$ as $n \to \infty$ is the Mandelbrot set for the exponential family $\lambda e^z$.  
We find this interesting because $Q_{n,c}$ is conjugate to $\Pnc$, but the Mandelbrot set for the exponential function is about as far away from the disk as a set can get:  for example, it contains the entire left half plane, plus infinitely long  ``hairs'' stretching roughly horizontally, toward infinity, in the direction of the positive real axis. 

Also, the Mandelbrot limit above follows from symmetry results on the sets $\MPnc$ found in \cite{Dierkfirst}.

Here, we obtain the result for $\MPnc$ by first investigating the Julia sets $\KPnc$ and $\JPnc$.  
Since $\MPnc$ tends to $\cDD$, the following result should make sense.

\begin{figure}
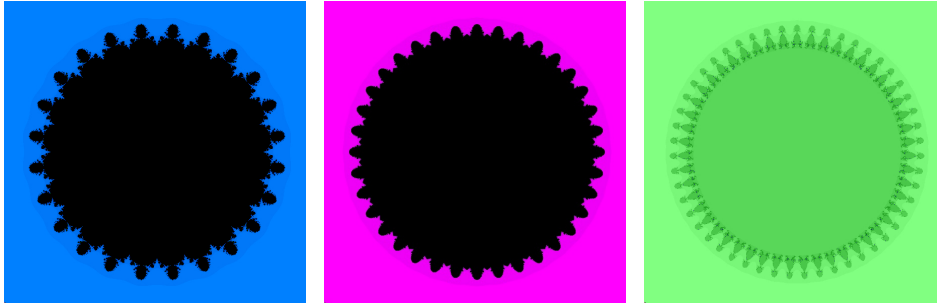

\drawfigMn
\drawfigJDn
\drawfigJSn
\caption{ Left: $\MPnc$ for $n=25$; Center: $\KPnc$ for $c=-.9, n=36$; 
Right: $\KPnc$ for  $n=51, c=1.10125i$.  
}
\end{figure}

\begin{thm} \label{thm:mainJPnc} Let $c \in \CC$.
\begin{enumerate}
\item If $c \in \CC\setminus\cDD$ then 
$$\lim_{n \to \infty} \JPnc =  \lim_{n\to\infty} \KPnc = S^1.$$
\item If $c \in \DD$, then
$$\lim_{n \to \infty} \JPnc = S^1  \  \text{ and } \ \lim_{n\to\infty} \KPnc = \cDD.$$
\item If $c \in S^1$, then if $\displaystyle\lim_{n\to\infty} \JPnc$ and/or $\displaystyle\lim_{n\to\infty} \KPnc$ (and/or any liminf or limsup) exists, it is contained in $\cDD$.

\hide{
\item If $c \in S^1$, then either
\begin{enumerate}
\item for all large $n$, $c \in \MPnc$, so  $\lim_{n\to\infty} \KPnc= \cDD$;
\item  for all large $n$, $c \notin \MPnc$,
so $\lim_{n\to\infty} \KPnc= S^1$; or
\item  neither of the above holds, so 
$\lim_{n\to\infty} \KPnc$ does not exist, but we can partition $\{ \KPnc \}$ into two subsequences according to whether $c \in \MPnc$ for each $n$.  One subsequence tends to ${\cDD}$, the other  tends to $S^1$.  Hence here we say
$$
\limsup_ {n\to\infty} \KPnc = {\cDD} \ 
\text{ and } \ 
\liminf_ {n\to\infty} \KPnc = S^1.
$$
\end{enumerate}
}

\end{enumerate}
\end{thm}

The case of $c \in S^1$ is the most delicate.  For any given $c \in S^1$, we may have $c \in M_n$ for some $n$ but not others.  We know the Julia sets may not extend outside $\cDD$ in the limit as $n$ grows, but for an individual $c \in S^1$, whether there is a limit inside the disk is uncertain.  This would be an interesting topic for future investigations.

For now, we move on to ask this same type of question for some other families of maps.
As a generalization of $\Pnc$, we are interested in the family
$$\Ftc(z) = z^t + c,$$
where $t \geq 2$ is real, and $z^t$ can be defined using, for example, the principal branch of the log: $\Ln(z)$, as in $z^t = e^{t \Ln(z)}$.
These functions are analytic, except on the negative real axis they are discontinuous.
We define a
filled Julia set $\KFtc$ as the closure of the set of points with bounded orbits,
the Julia set $\JFtc$ as its topological boundary, and the
Mandelbrot sets $\MFtc$ as the closure of the locus of $c$ with $0 \in \KFtc$.
Most of the nice results about polynomial/rational dynamics do not immediately apply, which is why to be safe we build ``closure'' into the definitions. (For example, we cannot say whether $\JFtc$ is the locus of nonnormality for the iterates of $\Ftc$, and we cannot say whether $\MFtc$ is the connectivity locus.)
However, with only minor changes to a few of the lemmas involved in proving Theorems~\ref{thm:mainMPnc} and \ref{thm:mainJPnc}, we see that the exact results of these theorems hold for $\Ftc$:

\begin{prop} \label{prop:mainFtc}
The limits in 
Theorem~\ref{thm:mainMPnc} and~\ref{thm:mainJPnc} may be taken as $t \to \infty$ through $\RR^+$ (rather than just as $n \to \infty$ through $\ZZ^+$).
\end{prop}

Finally, we find the most intriguing geometric limits when we investigate a family of rational maps which can be viewed as a peturbation of $\Pnc$: for $c,a \in \CC$ and $n \in \ZZ^+$, consider
$$\Rnca(z) = z^n + c + \frac{a}{z^n}.$$  

Since infinity is a superattracting fixed point, we still have a filled Julia set, $\KRnca$, defined as the set of points with bounded orbits, and the Julia set $\JRnca$ its topological boundary.   Note $0$ always maps to $\infty$ (for $a\neq 0$), thus the Julia sets lie in an annulus.
Though little is known about the dynamics for general $c,a$, we find that as $n$ grows, the Julia sets always tend to the unit circle:

\begin{thm} \label{thm:RcJ}  For any $c \in \CC$ and any $a \in \CC^*$,
$$ \lim_{n \to \infty} \JRnca = \lim_{n \to \infty} \KRnca = S^1.$$
\end{thm}

In investigating Mandelbrot sets, the case $c=0$ is most accessible, and has been most studied: Devaney's \cite{DevRatPert} contains a survey of known results regarding this family.  
For the case $c=0$, there happens to be only one free critical orbit. 
 Thus we may define the Mandelbrot set $\MRna$ as the set of all $a$ such that the free critical orbit does not escape to infinity.  
 %

For the case $c \neq 0$, there are two free critical orbits, thus the spectrum of dynamical behavior is quite rich.   One such case is investigated in \cite{DevRat}.   We let $\MRnca$ denote the set of $a$ such that for a fixed $n,c$, at least one free critical orbit does not escape to infinity.  These parameter sets are rife with mystery, but we did find that they settle into an easily describable limit as $n$ tends to $\infty$.

\begin{defn} \label{defn:limacon}
For any $c \geq 0$, define the set $L_{c}$ in $\CC$ by:
$$L_{c} = \left\{ \left( \frac{1}{4} + \frac{c}{2} e^{i\theta} +  \frac{1}{4}  e^{2i\theta} \right ) + \frac{(c^2 -1)}{4} \colon \theta \in [0, 2\pi) \right\}, $$  
and for any $c \in \CC$ define:
$$L_c =  e^{2 i \Arg{c}} L_{|c|}.$$
\end{defn}

Note $L_{c}$ is the polar lima\c con $r=(c + \cos\theta)/2$ shifted by $(c^2-1)/4$.

\begin{thm} \label{thm:MRnca}
For any $c \in \CC$, we have  $$\displaystyle \lim_{n\to \infty} \MRnca = L_c.$$ 
\end{thm}

This theorem includes the case $c=0$, and $L_0$ is the circle of radius $1/4$ centered at the origin, $S^1_{1/4}$.

\begin{figure}
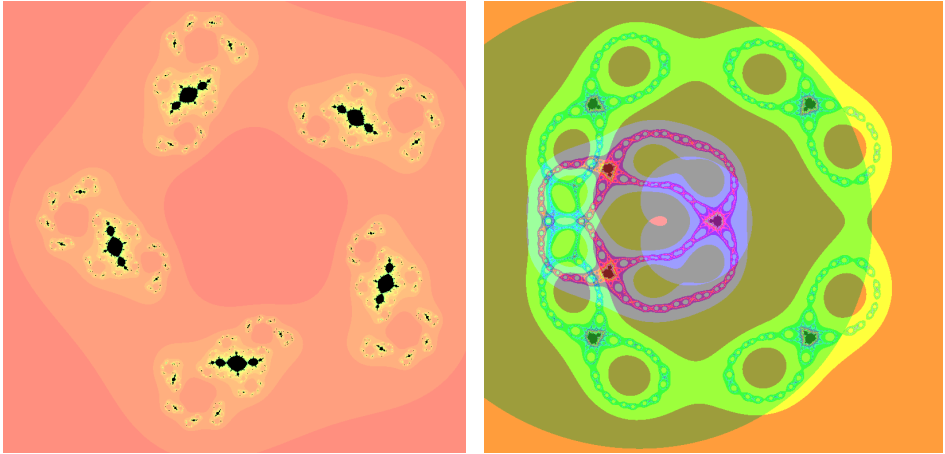

\drawfigJRnca
\drawfigMRnca
\caption{Left: $\JRnca$ for $n=5, c=1, a=0.0157+0.072i$. Right: $\MRnca$ for $n=8$, $c=0.25$.
}
\end{figure}

We found this avenue of study to be novel and exciting, and encourage the interested reader to study the geometric limit under degree growth of invariant sets of interest in other families of dynamical systems.

Proofs of these results, as well as more pictures, are located herein as follows: 
we  study   the polynomials $\Pnc$ in Section~\ref{sec:mainPnc}, the family $\Ftc$  in Section~\ref{sec:mainFtc}, and  the rational maps  $\Rnca$ in Section~\ref{sec:mainRnca}.  See Section~\ref{sec:background} for some definitions and notation used throughout the paper.

\subsection*{Acknowledgements} We thank Nicolet high school student Gabriel Ferns, for asking a question which led to this line of inquiry.  We also thank Hans Volker, Robert Devaney and Elizabeth Russell for helpful conversations, and Brian Boyd for the computer program detool which generated all of the Mandelbrot and Julia images in this paper, as well as inspired some of the results.

\section{Preliminaries}
\label{sec:background}

Here we explicitly describe what we mean when we say the limit of a sequence of sets exists.

Recall the \textbf{Hausdorff metric} $d_\HH(A,B)$ for two sets $A$ and $B$ of a metric space $(X,d)$ is 
\begin{align*}
d_\HH(A,B) & =  \max\{\sup_{a\in A} \inf_{b\in B} d(a,b), \  \sup_{b\in B} \inf_{a\in A}d(a,b)\} \\
 & =  \max\{ \sup_{a\in A} d(a,B),  \ \sup_{b\in B} d(b,A). \}
\end{align*}

Thus the Hausdorff distance between two sets $A, B$ can be thought of as the maximum of the maximum distance from $A$ to $B$ and the maximum distance from $B$ to $A$.    Consider the illustrative example $A=S^1$ and $B=\cDD$.  We have $d_\HH(S^1, \cDD) = 1$, since $0 \in \cDD$ and $d(0, S^1) = 1$, even though $d(e^{it}, \cDD) = 0$ for every $e^{it} \in S^1$.  Thus in each of our limit theorems, we need to supply two distance calculations.

If  $S_n$ is a sequence of compact subsets, then $S_n$ converges to compact subset $S$ and we write $\lim_{n \to \infty} S_n = S$
if
for all $\eps > 0$ there is an $N$ such that for all $n \geq N$ we have $d_\HH(S_n, S) < \eps$.

A compact metric space (like $\hat{\CC}$) endowed with the Hausdorff metric is complete (i.e., cauchy sequences converge). 

Mandelbrot and Julia sets for polynomials and rational functions are compact (\cite{Beardon}).
Since we define Mandelbrot and Julia sets for $\Ftc$ to be closed, we get they are compact if we can show they are bounded, which we will establish, at least for $n$ sufficiently large, in Section~\ref{sec:mainFtc}.
\medskip

Throughout the remainder of the paper we use the following notation:

{\bf Notation.}  For any $0 \leq r < R$, define the open annulus of inner radius $r$ and outer radius $R$ by:
$$ \An{r}{R} = \{ z \colon r < |z| < R \}.$$

{\bf Notation.} Denote the $\eps$-neighborhood of a set $S$ by $\mathcal{N}_\eps(S)$, and the ball about a point $z$ of radius $\eps$ by $B(z,\eps)$.

\section{Results for the family $\Pnc$}
\label{sec:mainPnc}

In this section, we first give a series of lemmas on the Julia sets for the maps $\Pnc(z) = z^n + c$, then prove Theorems~\ref{thm:mainMPnc} and~\ref{thm:mainJPnc}.  

We will distinguish the two cases (1) and (2) of Theorem~\ref{thm:mainJPnc} by the fundamental dichotomy of polynomial dynamics (see \cite{Beardon}):  

{\bf The Fundamental Dichotomy.} \textit{Since $\Pnc$ is unicritical, for any given $n$ and $c$ either 
\begin{enumerate}
  \item[(a)] the orbit of $0$ escapes, in which case $\KPnc = \JPnc$ is a Cantor set, or 
  \item[(b)] the orbit of $0$ is bounded, in which case $\KPnc$ is connected.  
  \end{enumerate}
}

Case (1): Suppose $\abs{c} > 1$. We will show that for sufficiently large $n$, $0 \notin \KPnc$, hence (a) holds and $c \notin M_n$.

Case (2): On the other hand, suppose $\abs{c}<1$. We will show that for large $n$, $0 \in \KPnc$, hence (b) holds and $c \in M_n$.

\medskip

\hide{
Along the way toward establishing Theorem~\ref{thm:mainMPnc}, we will actually find:
 
\begin{prop} \label{prop:MPncconverge}
The exterior of $\MPnc$ (the Cantor locus) tends to $\CC \setminus \DD$ as $n\to \infty$, and the central hyperbolic component of $\MPnc$ (the set of $c$ such that $\Pnc$ has an attracting fixed point) tends to $\DD$.  
\end{prop}
}

\subsection{Julia sets}
We first establish one half of the distance calculation for each case of $|c|$: showing the distances from every point of the Julia sets and filled Julia sets to $S^1$ or $\cDD$, as appropriate, tend to $0$ as $n \to \infty$.

We begin by showing that regardless of $c$, the immediate basin of infinity of $\Pnc$ encompasses the entire exterior of the unit disk as $n \to \infty$:

\begin{lem} \label{lem:KinDgeneral}
For any $c \in \CC$, and any $\eps >0$, there is an $N \geq 2$ such that for all $n \geq N$,
we have $K(\Pnc) \subset \DD_{1+\eps}$.
\end{lem}

\begin{proof}
Let $z \in \CC \setminus \cDD_{1+\eps}$.  Choose $B > \max\{1, |c|\}$.  
Choose $N \geq 2$ so that $|z|^N > 2B$ and $B^{N-2} > 2$.  Let $n \geq N$.  We claim 
$|\Pnc^m(z) | \geq B^m$ for all $m \geq 1$.  

First, $$|\Pnc(z)| = |z^n+c|  \geq |z|^n-|c| > 2B - B = B.$$

Now suppose for some $m \geq 1$, we know $|\Pnc^m(z) | \geq B^m$.  Then
\begin{align*}
|\Pnc^{m+1}(z)|
&=|(\Pnc^m(z))^n+c|\\
&\ge|\Pnc^m(z)|^n-|c|\\
&>B^{mn}-B\\
&>B^{mn-m+1}-B\\
&=B^{m+1}(B^{m(n-2)}-B^{-m})\\
&>B^{m+1}(2^m-B^{-m})\\
&>B^{m+1}(2^m-1)\\
&\ge B^{m+1}.
\end{align*}
So, by induction, $|\Pnc^m(z)|>B^m$ for all $m \geq 1$.  Since $B > 1$, the orbit of $z$ under $\Pnc$ escapes to infinity, thus $z \notin K(\Pnc)$.  Thus, $K(\Pnc) \subset \DD_{1+\eps}$.

\end{proof}

Note this lemma immediately yields (3) of Theorem~\ref{thm:mainJPnc}. 

\begin{figure}
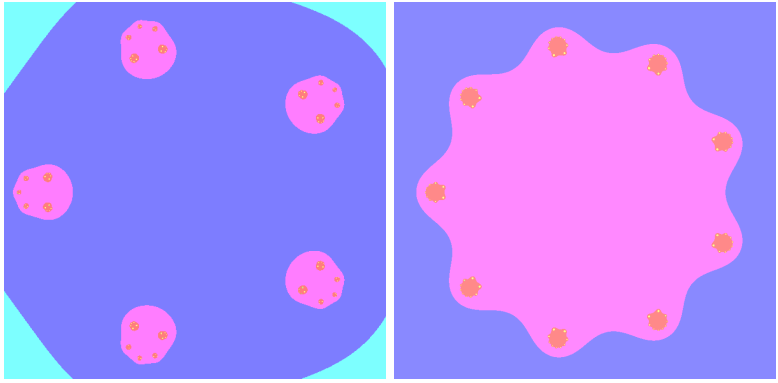

\drawfigKPoutnlow
\drawfigKPoutnhi
\caption{\label{fig:JPncout}
$\KPnc$ for $c=2$, $n=5$ (L) and $n=9$ (R).}
\end{figure}

Now we examine the case $\abs{c} >1$.    Figure~\ref{fig:JPncout} displays some Julia sets with $\abs{c}>1$ and $n$ large enough that these $c \notin M_n$, in fact $\JPnc = \KPnc$ and these are Cantor Julia sets, contained in an annulus which tends to $S^1$ as $n \to \infty$.

\begin{lem} \label{lem:KinDgt}
Let $c\in\CC\setminus\cDD$ and set $\eps=1-|c|$. Then, for every $\eta \in (0,\eps]$, there is an $N\ge 2$ where $\DD_{1-\eta/2}\subset\CC\setminus K(\Pnc)$ for all $n\ge N$.
\end{lem}

\begin{proof}
Choose $N_1$ by Lemma~\ref{lem:KinDgeneral} so large that for all $n \geq N_1$ we have $\KPnc \subset \DD_{1+\eta/2}.$
Now let $|z|<1-\eta/2$, and chose $N \geq N_1$ so that $|z|^N<\eta/2$. Then for all $n \geq N$,
\begin{align*}
|\Pnc(z)|&=|z^n+c|\\
&\ge|c|-|z|^n\\
&=1+\eps-|z|^n\\
&\ge 1+\eta-|z|^n\\
&>1+\eta-\eta/2\\
&=1+\eta/2.
\end{align*}
By Lemma~\ref{lem:KinDgeneral}, $\Pnc(z)\notin K(\Pnc)$ so neither is $z$.
\end{proof}

Combining Lemmas~\ref{lem:KinDgeneral} and~\ref{lem:KinDgt}, we obtain:

\begin{cor} \label{cor:KoutDinA}
Let $c \in \CC \setminus \cDD$ and set $\eps = |c| -1$. Then, for every $\eta \in (0, \eps]$, there is an $N \geq 2$ such that for all $n \geq N$, we have $\JPnc = \KPnc \subset \An{1-\eta/2}{1+\eta/2}$.
\end{cor}

This corollary provides half of the needed distance calculation for case (1) of Theorem~\ref{thm:mainJPnc}:  $\displaystyle\inf_{z \in \JPnc = \KPnc} d(z, S^1)$ tends to $0$ as $n\to\infty$.

\medskip

\begin{figure}
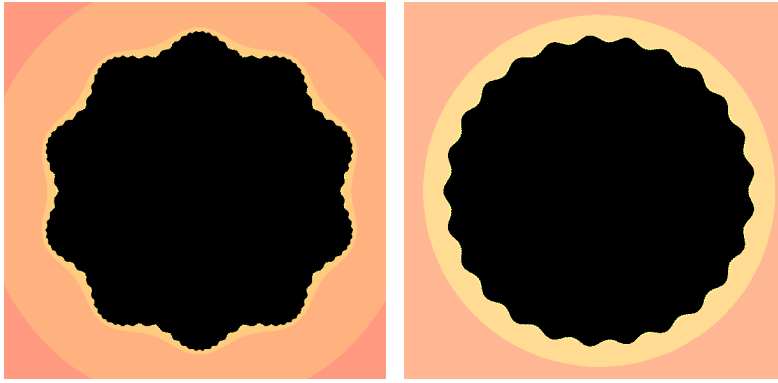

\drawfigKPinnlow
\drawfigKPinnhi
\caption{\label{fig:JPncinD}
$\KPnc$ for $c=0.5$, $n=10$ (L) and $n=25$ (R).}
\end{figure}

Next we consider the case $\abs{c} < 1$.  Figure~\ref{fig:JPncinD} shows some Julia sets for $c \in \DD$ and $n$ large enough that $c \in M_n$. In fact here $c \in M_n^{\circ}$ and these filled Julia sets illustrate how $\KPnc$ can tend to $\cDD$ as $n\to\infty$ while $\JPnc$ tends to $S^1$. 

\begin{lem} \label{lem:DinK}
Let $c\in\DD$, set $\eps = 1 - |c|$.  For any $\eta \in (0, \eps]$,  there is an $N\ge 2$ so that $\DD_{1 - \eta/2}\subset K(\Pnc)^{\circ}$ (the interior of $\KPnc$, i.e., $\KPnc \setminus \JPnc$) for all $n\ge N$.
\end{lem}

\begin{proof}
Let $z\in\DD_{1-\eta/2}$.  Since $\eta > 0$, there is an $N\ge 2$ where $(1-\eta/2)^n<\eta/2$ for all $n\ge N$.  For such $n$,

\begin{align*}
|\Pnc(z)|&\le|z|^n+|c|\\
&<\frac\eta2+1-\eps \leq \frac{\eta}{2} + 1 - \eta \\
&=1-\frac\eta2.
\end{align*}

So $\Pnc(\DD_{1-\eta/2})\subsetneq \DD_{1-\eta/2}$, and thus $\DD_{1-\eta/2}\subset K(\Pnc)^{\circ}$.
\end{proof}

Now combining Lemmas~\ref{lem:KinDgeneral} and~\ref{lem:DinK}, we obtain:

\begin{cor} \label{cor:JinDinA}
Let $c \in \DD$ and set $\eps = 1 - |c|$. For any $\eta \in (0, \eps]$, there is an $N \geq 2$ such that for all $n \geq N$, we have $\JPnc \subset \An{1-\eta/2}{1+\eta/2}$
and $\DD_{1-\eta/2} \subset \KPnc \subset \DD_{1+\eta/2}$.
\end{cor}

This corollary provides half the distance calculation needed for case (2) of Theorem~\ref{thm:mainJPnc}:  $\displaystyle\inf_{z \in \JPnc} d(z, S^1)$ tends to $0$ as $n\to\infty$, and  $\displaystyle\inf_{z \in \KPnc} d(z, \DD)$ tends to $0$ as $n \to \infty$.


We turn to supplying the other halves of the distance calculation needed for cases (1) and (2).  Note that since $\JPnc \subset \KPnc$, we need only show that for all $s \in S^1$, 
the distance $d(s, \JPnc)$ tends to $0$ as $n \to \infty$.

\hide{
\begin{lem} \label{lem:JinA}
Let $c\in\CC$.  Then for any $\eps>0$, there is an $N\ge2$ such that $J(P_{n,c})\subset\DD_{1+\eps/2}\setminus\DD_{1-\eps/2}$ for all $n\ge N$.
\end{lem}
}  

We do so by taking advantage of the $n$-fold symmetry of $\JPnc$:

\begin{lem} \label{lem:rtpi}
Let $z\in J(P_{n,c})$.  If $\omega$ is a $n^{th}$ root of unity, then $\omega z\in J(P_{n,c})$.
\end{lem}

\begin{proof}
Let $\omega$ be an $n^{th}$ root of unity, and suppose $z\in J(P_{n,c})$.  Then
$$
P_{n,c}(\omega z)=(\omega z)^n+c =\omega^nz^n+c =z^n+c.
$$
Since $\JPnc$ is totally invariant (see \cite{Beardon}), $\omega z$ must be in $\JPnc$.
\end{proof}


\begin{lem} \label{lem:JdinS}
Let $c\in\CC \setminus S^1$, set $\eps = \abs{1 - |c|}$ and let $\eta \in (0, \eps]$. 
Then there is an $N\ge2$ such that for all $n\ge N$ and for any $e^{i \theta} \in S^1$, $B(e^{i \theta},\eta)\cap J(P_{n,c})\ne\emptyset$.
\end{lem}

\begin{proof}

By Corollaries~\ref{cor:KoutDinA} and~\ref{cor:JinDinA}, there is an $N_1\ge 2$ where  $J(P_{n,c})\subset\An{1-\eta/2}{1+\eta/2}$ for all $n\ge N_1$.

Let $e^{i\theta} \in S^1$ and let $\alpha>0$ be the angle so that $U=\{re^{i\tau}| r>0, \theta-\alpha<\tau<\theta+\alpha\}\cap \An{1-\eta/2}{1+\eta/2}$ is contained in $B(e^{i\theta},\eta)$. Note the same $\alpha$ works for each different $\theta$.

Let $\omega_n = e^{2\pi i/n}$ for any $n$.  Now, choose $N\ge N_1$ so that 
$2\pi/N < \alpha$, hence
$2\pi/n < \alpha$ for all $n \geq N$. Again, note $N$ is independent of $\theta$, i.e., the same $N$ works for every $\theta$.

Since for any $n$, $\JPnc$ is nonempty (\cite{Beardon}), choose $z_{n}\in J(P_{n,c})$ for each $n \geq N$.  Then for some $j_n \in \{ 1, \ldots, n-1\}$ we have $\omega_{n}^{j_n} z_n \in U\subset B(e^{i\theta},\eta)$.

So for all $n\ge N$, $B(e^{i\theta},\eta)\cap J(P_{n,c})\ne\emptyset$.
\end{proof}

\begin{proof}[Proof of Theorem~\ref{thm:mainJPnc}]

(\textit{Case (1)}) Let $c\in\CC\setminus\overline\DD$.  Since by Corollary~\ref{cor:KoutDinA}, $\JPnc = \KPnc$ for $n$ sufficiently large, we write only $\JPnc$ below. Let $\eps>0$ and assume $\eps < |c|-1$.  By Corollary~\ref{cor:KoutDinA},  there is an $N_1$ such that for all $n \geq N_1$, $\JPnc \subset\An{1-\eps/2}{1+\eps/2}$, and thus for all $z \in \JPnc$, 
$$
d(z, S^1) =  \inf_{s\in S^1}|z-s|<\eps/2.
$$
By Lemma~\ref{lem:JdinS}, there is an $N_2 \geq N_1$ so that for all $n \geq N_2$ and for all $s \in S^1$,
$$
d(s, \JPnc) =  \inf_{z\in\JPnc}|z-s|<\eps.
$$
So we have
\begin{align*}
d_\mathcal H(\JPnc,S^1)&=\max\bigg\{\sup_{z\in\JPnc}d(z,S^1),\sup_{s\in S^1}d(s,\JPnc)\bigg\}\\
                      &<\max\{\eps/2,\eps\}\\
                      &=\eps.
\end{align*}
Hence we have $\lim_{n \to \infty} \JPnc = \lim_{n\to\infty} \KPnc = S^1$.

(\textit{Case (2)}) Suppose $c \in \DD$.  Let $\eps > 0$ and assume $\eps < 1-|c|$.  By Corollary~\ref{cor:JinDinA}, there is an $N_1$ such that for all $n \geq N_1$, we have 
$\JPnc \subset \An{1-\eps/2}{1+\eps/2}$ and  $\DD_{1-\eps/2} \subset \KPnc \subset \DD_{1+\eps/2}$.  We can establish the limit for $\KPnc$ using only the preceeding statement.  First, note
for all $w \in \KPnc$,
$$
d(w, \cDD) = \inf_{y \in \cDD} |w-y| < \eps/2.
$$
Now, consider $y \in \cDD$.  If $y \in \DD_{1-\eps/2}$, then $y \in \KPnc$, hence  $d(y, \KPnc) = 0$.
If $1 - \eps/2 \leq |y| \leq 1$, then $d(y, \DD_{1-\eps/2}) \leq \eps/2$, and since this disk is contained in $\KPnc$, we have $d(y, \KPnc) \leq \eps/2$.  Thus for any $y \in \cDD$, we have
$$
d(y, \KPnc) = \inf_{w \in \KPnc} |w-y| \leq \eps/2 < \eps.
$$
Thus, $d_\HH (\KPnc, \cDD) < \eps$, and hence $\lim_{n \to \infty} \KPnc = \cDD$.

Now we turn to $\JPnc$.  By choice of $N_1$ above, we know for all $z \in \JPnc$,
$$
d(z, S^1) =  \inf_{s\in S^1} |z-s| < \eps/2.
$$

Now, by Lemma~\ref{lem:JdinS}, there is an $N_2 > N_1$ such that for all $n \geq N_2$, 
for any point $s=e^{i \theta} \in S^1$, 
$$
d(s, \JPnc) = \inf_{z\in\JPnc} |z-s| < \eps/2.
$$

Thus 
$d_\HH (\JPnc, S^1) < \eps/2$, and 
hence $\lim_{n \to \infty} \JPnc = S^1$.

Case (3) follows directly from Lemma~\ref{lem:KinDgeneral}.

\end{proof}

\subsection{Mandelbrot limits}
We now turn to the Mandelbrot limits, establishing
Theorem~\ref{thm:mainMPnc}.

\begin{figure}
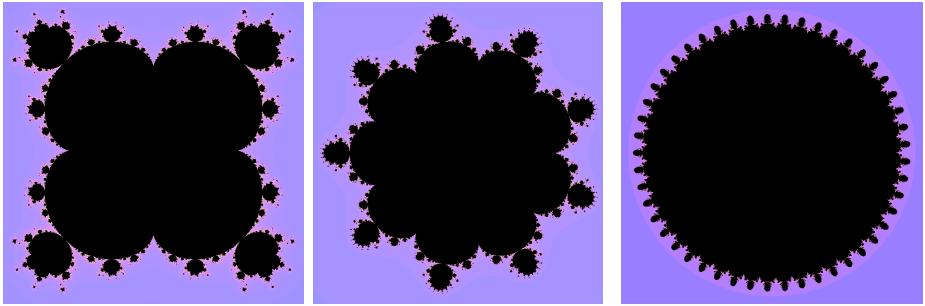

\drawfigMPnlo
\drawfigMPnmed
\drawfigMPnhi
\caption{\label{fig:MPnc}
$\MPnc$ for $n=5, 10, $ and $50$ (L to R).}
\end{figure}

Figure~\ref{fig:MPnc} illustrates how $\MPnc$ seems to tend to the closed unit disk as $n$ grows, and suggests that the central hyperbolic component grows into the entire unit disk as $n$ grows (which Lemma~\ref{lem:DinK} supports).

The needed distance calculations for the main Theorem will follow pretty easily from the lemmas on Julia sets we proved above.

\begin{cor} \label{cor:DinM}
For every $\eps\in(0,1)$, there is an $N\ge2$ so that $\overline{\DD_{1-\eps}}\subset \MPnc$ for every $n\ge N$.
\end{cor}

\begin{proof}

First, note that $c=0 \in \MPnc$ for every $n \geq 2$.

Now, given $\eps \in (0,1)$, let $c \in \overline{\DD^*_{1-\eps}}$.  Then $|c| \leq 1-\eps$.  Set $\eps' = 1-|c|$, so $0 < \eps' \leq \eps$. Then by Lemma~\ref{lem:DinK}, with $\eta = \eps'$, 
there is an $N \geq 2$ such that for all $n \geq N,$ we have $0 \in \DD_{1-\eps'/2} \subset K(\Pnc)$.  Hence $c \in \MPnc$.

\end{proof}

\begin{cor} \label{cor:MinD}
For every $\eps>0$, there is an $N\ge2$ so that $\MPnc\subset \overline{\DD_{1+\eps}}$ for every $n\ge N$.
\end{cor}

\begin{proof}
Let $\eps > 0$.  Suppose $|c| > 1 + \eps$.   Then by Lemma~\ref{lem:KinDgeneral}, there is an $N \geq 2$ such that for all $n \geq N$, $K(\Pnc) \subset \DD_{1+\eps}$ but $\Pnc(0) = c \notin \DD_{1+\eps}$, hence $c \notin K(\Pnc)$, hence $c \notin \MPnc$.
\end{proof}

\hide{
\begin{proof}[Proof of Proposition~\ref{prop:MPncconverge}]
The first statement follows from Lemma~\ref{lem:KinDgeneral}. 

Lemma~\ref{lem:DinK} shows that for every $c \in \DD$, for  $\eps>0$ sufficiently small, $\Pnc$ maps $\DD_{1-\eps}$ strictly inside itself, which implies that in that $\Pnc$ has an attracting fixed point in that disk (\cite{Beardon}).

Thus, we have that the main hyperbolic component (the one containing $c=0$) tends to $\DD$.

\end{proof}
} 


\begin{proof}[Proof of Theorem~\ref{thm:mainMPnc}]
Let $\eps>0$.  By Corollaries~\ref{cor:DinM} and \ref{cor:MinD}, there are $N^+$ and $N^-$ where $\MPnc\subset{\cDD_{1+\frac{\eps}{2}}}$ for all $n\ge N^+$, and ${\cDD_{1-\frac{\eps}{2}}}\subset \MPnc$ for all $n\ge N^-$.  Set $N=\max\{N^+,N^-\}$.  Then for all $n \geq N$,
for any $c \in \MPnc$, we have
$$
d(c, \cDD) = \inf_{y \in \cDD} |c-y| < \eps/2.
$$

Now consider $y \in \cDD$.  If $|y| \leq 1-\eps/2$, then $y \in \MPnc$, hence
$d(y, \MPnc) = 0$.  If $1 - \eps/2 < |y| \leq 1$, then $d(y, \cDD_{1-\eps/2}) < \eps/2$, and since $\cDD_{1-\eps/2} \subset \MPnc$, we have $d(y, \MPnc) < \eps/2$.  Hence for any $y \in \cDD$, we have 
$$
d(y, \MPnc) = \inf_{c \in \MPnc} |c-y| < \eps/2.
$$

Thus, $d_\HH (\cDD, \MPnc) < \eps/2$, hence
 $\lim_{n\to\infty} \MPnc = \cDD$.

\end{proof}

We also provide an alternative way to say that $\MPnc$ tends to $\cDD$:

\begin{lem} \label{lem:dnsS}
There is a dense set $Q\subset S^1$ where for any $N\ge 2$,
\[
Q\subset \( \bigcup_{n\ge N}\MPnc \) =: A_N
\]
\end{lem}

\begin{proof}
Let $n\ge 2$.  Here, the solution set of 
$c^n+c=0$ is $Q_n\cup\{0\}$ where
\[
Q_n=\{z|z\text{ is an }(n-1)\text{th root of }-1\}.
\]
If $c\in Q_n$, then the orbit of $\Pnc(0)$ is $\{0,c,0,c,\dots\}$ and thus bounded, so $c\in \MPnc$.
See Figure~\ref{fig:PncPertwo}.

\begin{figure}
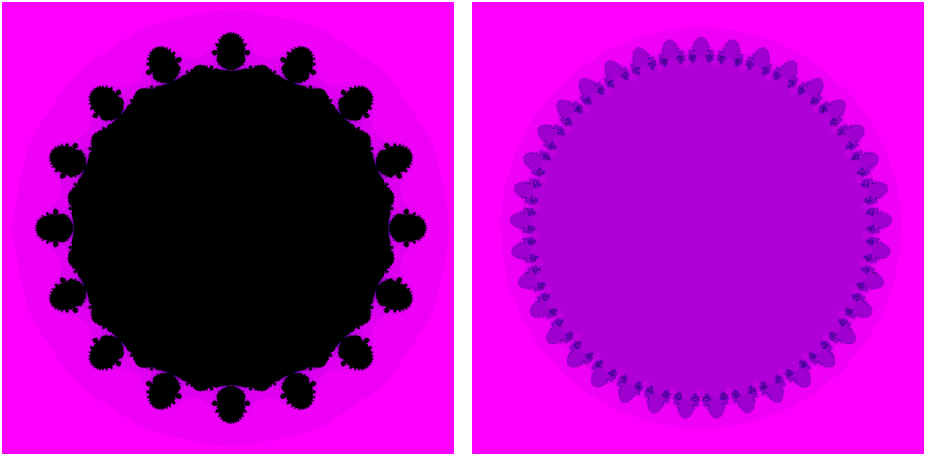

\drawfigKPincone
\drawfigKPoutcone
\caption{\label{fig:PncPertwo}
$\KPnc$ for $n=16, c=-0.9$ (L) and $c=-0.9i, n=37$ (R). 
For many $|c| \approx 1$, it appears that  $c \in \MPnc$ for some $n$, but not for others.
}
\end{figure}

Define $Q = \{ e^{2\pi qi} \colon q\in\QQ\cap(0,1) \}$. Let $c \in Q$.
  Then there is an $n$ such that for all $k \geq 1,  c \in Q_{k(n-1)+1} \subset M_{k(n-1)+1}$.  Thus $c\in A_N$, for all $N \geq 2$.  So $Q\subset A_N$ for all $N\ge 2$, and clearly $Q$ is dense in $S^1$.
\end{proof}

\begin{thm}
The accumulation set of the sequence $\{ \MPnc \}_{n \to \infty}$ is the closed unit disk, i.e.,
\[
\cDD = \( \bigcap_{N\ge2} \overline{\bigcup_{n\ge N} \MPnc} \) =: A.
\]
\end{thm}

\begin{proof}
\begin{list}{\label itemi}{\leftmargin=1em}
\item[$(\cDD\subset A)$] Let $z\in\DD$.  By Corollary~\ref{cor:DinM}, there is an $N$ so that $z\in \MPnc$ for all $n\ge N$, thus $z\in A$.  If $z\in\partial\DD=S^1$, then by Lemma~\ref{lem:dnsS}, since $Q$ is dense in $S^1$, $S^1\subset\overline{A_N}$ for all $N\ge 2$, and since $A=\bigcap_{N\ge 2}\overline{A_N}$, $z\in A$ and $\cDD\subset A$.\\
\item[$(A\subset\cDD)$] Let $z\in\CC\setminus\cDD$.  By Corollary~\ref{cor:MinD}, there is an $N\ge 2$ where $z\notin \MPnc$ for all $n\ge N$.  Thus $z\notin A$.  So contrapositively, $A\subset\cDD$.
\end{list}
\end{proof}

\section{Results for the family $\Ftc$}
\label{sec:mainFtc}

For one generalization of $\Pnc$, we are interested in the family
$$\Ftc(z) = z^t + c,$$
where $t > 1$ is real, and $z^t$ can be defined using, for example, the principal branch of the log: $\Ln(z)$, as in $z^t = e^{t \Ln(z)}$.

These functions are analytic, except on the negative real axis they are discontinuous.
We may still define the
{filled Julia set} as 
$$\KFtc = \overline{\{ z \colon \{\Ftc^m(z)\}_{m=1}^{\infty}  \text{ is bounded}\} },$$
the {Julia set} as the topological boundary, $\JFtc= \partial \KFtc$, and the
{Mandelbrot sets} as 
$$\MFtc= \overline{\{ c  \colon \{\Ftc^m(0)\}_{m=1}^{\infty} \text{ is bounded} \}}.$$

\begin{figure}
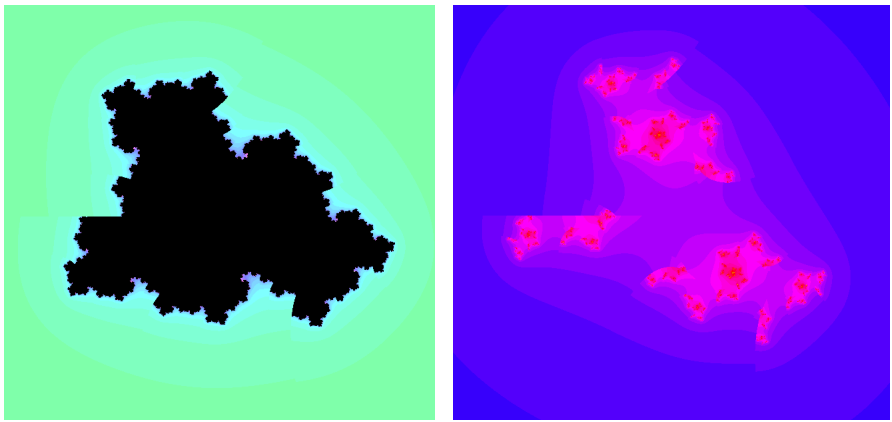

\drawfigKFtcin
\drawfigKFtcout
\caption{\label{fig:JFtc} $\KFtc$ for $t=2.5$, $c=-0.5+0.5i$ (L) and $c=i$ (R).}
\end{figure}

\begin{figure}
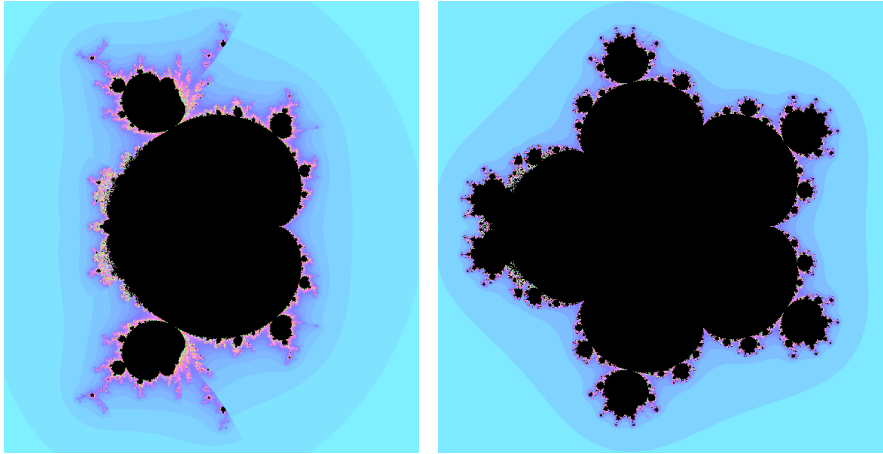

\drawfigMFtcnlo
\drawfigMFtcnhi
\caption{\label{fig:MFtc} $\MFtc$ for $n=2.5$ (L) and $6.25$ (R)}
\end{figure}

This discontinuity is reflected in the Julia sets, see Figure~\ref{fig:JFtc} and Mandelbrot sets, see Figure~\ref{fig:MFtc}.

Irregardless, a brief inspection shows that proofs of all of the results on $\Pnc$ of Section~\ref{sec:mainPnc} easily generalize to $\Ftc$.
One need only change $n, N \in \NN$ to $t, T \in \RR$ and for any integer $p$, replace $|z|^p$ by $|z^{\rho}|$ for the appropriate real number $\rho$.  We use the observation that even though $z \mapsto z^t$ is not continuous, the function $z \mapsto |z^t|$ is continuous (and not multiply defined).

The one exception is that Lemma~\ref{lem:dnsS} needs no change, we can apply it as is. In adapting the proofs, the references to $A_N$ are unchanged, but other integers and lemmas are changed to their $t$-versions.  All the other results hold with the changes mentioned above.  The Julia set and Mandelbrot set bounds for large $n$ imply the Julia sets are bounded, and since we defined them to be closed, they are compact.  Thus we obtain Proposition~\ref{prop:mainFtc}.

For brevity sake, we leave the verification of this statement to the reader.

\subsection*{Future Research} Images suggest that $\MFtc$ varies continuously with $t$ for all $t >1$ (or possibly over a larger range of $t$ values), and that under some conditions on $c$, the Julia sets also vary continuously with $t$.

\section{Results for the family $\Rnca$}
\label{sec:mainRnca}

Now we extend the results for the polynomial family $\Pnc(z) = z^n+c$ to the rational family 
$$\Rnca(z) = z^n + c + \frac{a}{z^n}.$$ 
Devaney's \cite{DevRatPert} contains a survey of known results regarding this family.
Little is known when $c \neq 0$ (except in the case $c$ is periodic and $a$ is small).
Since infinity is a superattracting fixed point, we still have a filled Julia set, $\KRnca$, defined as the set of points with bounded orbits, and the Julia set $\JRnca$ is its topological boundary. Both $K$ and $J$ are compact and nonempty (\cite{Beardon}).

These maps have more critical points than $\Pnc$:  infinity is critical and fixed, zero is critical and maps to infinity, and there are also critical points at the $(2n)^{th}$ roots of $a$:  $a^{1/2n}$. 
Note that the $2n$ critical points $a^{1/2n}$ only generate two critical values: $v^{\pm} = c \pm 2 \sqrt{a}$.

\subsection{Julia sets}

\begin{figure}
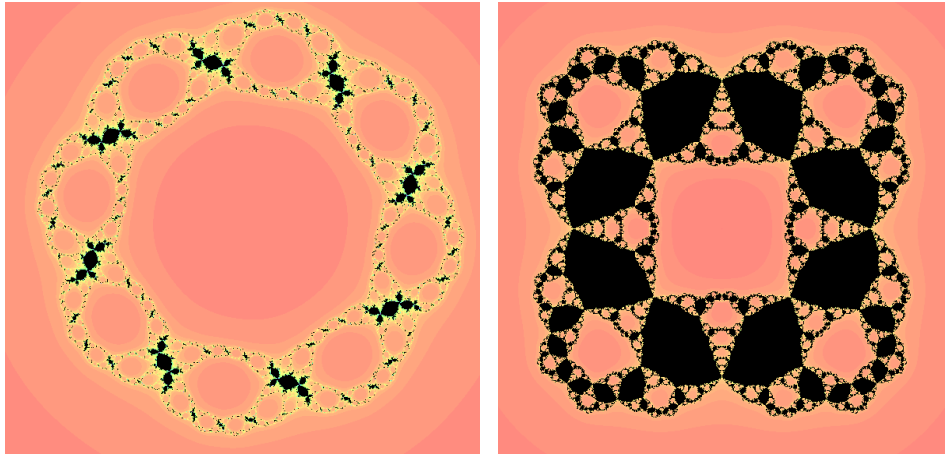

\drawfigKRnaexmp
\drawfigKRncaexmp
\caption{\label{fig:JRnca} Left: $\KRnca$ for $n=4, c=0, a=0.0422-0.147i$. 
Right: $\KRnca$ for $n=4, c=.5, a=-0.01253761$.}
\end{figure}

We start with lemmas restricting the location of the Julia sets to an annulus tending to $S^1$ as $n \to \infty$.  See Figure~\ref{fig:JRnca}.

\begin{lem} \label{lem:RcKinD}
For any $c \in \CC$ and any $a \in \CC$, given any $\eps > 0$, there is an $N \geq 2$ such that for all $n \geq N$, we have $\KRnca \subset \DD_{1 + \eps}$.
\end{lem}

\begin{proof}
Let $z \in \CC \setminus \cDD_{1+\eps}$.
Let $B > \max\{1, |a|, |c|\}$.  Choose $N \geq 2$ such that $|z|^N > 3B$ and $B^{n-2} > 2^2$.  Let $n \geq N$.  We claim that $|\Rnca^m(z)| \geq B^m$ for all $m \geq 1$.
First, 
\begin{align*}
|\Rnca(z)| &= \left | z^n + c + \frac{a}{z^n}\right | 
 \geq |z|^n - |c| - \frac{|a|}{|z|^n} \\
& > 3B - B - \frac{|a|}{2B} 
 = 2B - \frac{|a|}{B}
= B + \left (B - \frac{|a|}{B} \right ) > B.
\end{align*}

Now, suppose for some $m \geq 1$, $|\Rnca^m(z)| \geq B^m$.  Then
\begin{align*}
|\Rnca^m(z)| & = \left | (\Rnca^m(z))^n + c + \frac{a}{(\Rnca^m(z))^n} \right | \\
& \geq   | \Rnca^m(z)|^n - |c| - \frac{|a|}{ | \Rnca^m(z)|^n}\\
& > B^{mn} - B - \frac{|a|}{B^{mn}}\\
& > B^{mn-m+1} - B - \frac{|a|}{B^{mn}}\\
& = B^{m+1} \left ( B^{m(n-2)} - \frac{|a|}{B^{mn+m+1}} - B^{-m} \right )\\
& \geq B^{m+1} (2^{2m}-1-1) \\ & \geq B^{m+1}.
\end{align*}

Thus by induction, $|\Rnca^m(z)| \geq B^m$ for all $m \geq 1$.  Since $B>1$, 
the orbit of $z$ under $\Rnca$ escapes to infinity, hence $z \notin \KRnca$.
\end{proof}

\begin{lem} \label{lem:RcDoutK}
For any $c \in \CC$ and any $a \in \CC^*$, given any $\eps > 0$ there is an $N \geq 2$ such that for all $n \geq N$, we have $\DD_{1-\eps} \subset \CC \setminus \KRnca$.
\end{lem}
 
 \begin{proof}
Let $z \in \DD_{1-\eps}$.  Let $N_1$ be given by Lemma~\ref{lem:RcKinD} for this $a,c, \eps$, so that for all $n \geq N_1$ we have $\KRnca \subset \DD_{1+\eps}$.  Now, since $|z| < 1$, choose $N \geq N_1$ such that $|z|^N \leq \abs{a}/(2+\eps+\abs{c})$.  Then for any $n \geq N$, 

\begin{align*}
|\Rnca(z)| & \geq \frac{|a|}{|z|^n} - \abs{c} - \abs{z}^n \\
& \geq  \frac{|a|}{|z|^n} - \abs{c} - 1\\
& \geq 1 + \eps.
\end{align*}

Hence $\Rnca(z) \notin \KRnca$, thus $z \notin \KRnca$. 
 \end{proof}
 
Combining these lemmas we obtain:
 
 \begin{cor} \label{cor:RcKinA}
For any $c \in \CC$ and any $a \in \CC^*$, given any $\eps > 0$, there is an $N \geq 2$ such that for all $n \geq N$, we have $\KRnca \subset \An{1-\eps}{1+\eps}$. 
 \end{cor}

  
Since $J \subset K$,  we now can conclude that (if $c \in \CC$ and $a \in \CC^*$), both $\inf_{z \in \JRnca} d(z, S^1)$ and $\inf_{z \in \KRnca} d(z, S^1)$ tend $0$ as $n\to\infty$, which is half the distance calculation needed for Theorem~\ref{thm:RcJ}.  
We now establish the other distance, taking advantage of symmetry in $\Rnca$:

 \begin{lem} \label{lem:RcSinJ}
 For any $c \in \CC$ and any $a \in \CC$,  if $\alpha$ is a primitive $n^{th}$ root of unity, then
 for any point $z \in \JRnca$ we have $\{ \alpha^j z \colon j \in \{ 0, 1, \ldots, n -1\} \}$ are $n$ distinct points in $\JRnca$, spread evenly around the circle of radius $\abs{z}$.
\end{lem}

\begin{proof}
Since $\alpha$ is a primitive $n^{th}$ root of unity, $\alpha^n = 1$ but $\alpha^j \neq 1$ for any $j < n$, so since $0 \notin \JRnca$ implies $z \neq 0$, we have the set of $\alpha^j z$ a set of $n$ distinct points, and note for each $j$, $(\alpha^j)^n = 1.$  Thus, observe
 $$\Rnca( \alpha^j z) = (\alpha^j z)^n + c + \frac{a}{(\alpha^j z)^n} = (\alpha^j)^n z^n + c+ \frac{a}{(\alpha^j)^n z^n} = \Rnca(z).$$  
 Hence each $\alpha^j z$ maps to the image of $z$.  Thus since $J$ is totally invariant (\cite{Beardon}), if $z$ is in $\JRnca$ then so is each $\alpha^j$, and each has norm $\abs{z}$.
\end{proof}

\begin{lem} \label{lem:S1inJ}
Let $c \in \CC, a \in \CC$.  Let $\eps > 0$.  Then there is an $N \geq 2$ such that for all $n \geq N$, for any $s = e^{i\theta}$ we have $B(s, \eps) \cap \JRnca \neq \emptyset$.
\end{lem}

\begin{proof}
Now, we know $\JRnca$ is nonempty, since $\Rnca$ is a rational map (\cite{Beardon}).  Hence there exists a point $z_{n} \in \JRnca$ for any $n \geq 2$.  
Since Lemma~\ref{lem:RcSinJ} gives us that $z_{n} \in \JRnca$ implies the points $\alpha^j z_{n}$ are also in $\JRnca$, we have $n$ evenly distributed points along the circle of radius $\abs{z_{n}}$ which are in $\JRnca$ and 
by Corollary~\ref{cor:RcKinA}, as $n \to \infty$ we have $\abs{z_{n}} \to 1$. 

The same argument as in the proof of Lemma~\ref{lem:JdinS} will finish the proof.
\end{proof}
 
 \begin{proof}[Proof of Theorem~\ref{thm:RcJ}]
First recall for any $n, c,a$ we have $\JRna \subseteq \KRna$.  

In direct analogy to the proof of Theorem~\ref{thm:mainJPnc}, Corollary~\ref{cor:RcKinA} implies that both $\displaystyle \inf_{z \in \JRnca} d(z, S^1)$ and $\displaystyle \inf_{z \in \KRnca} d(z, S^1)$ tend $0$ as $n\to\infty$.

Lemma~\ref{lem:S1inJ} implies that 
for any $s \in S^1$, we get $d(s, \JRnca)\to 0$ as $n \to \infty$, and since $\JRnca \subset \KRnca$, this implies $d(s, \KRnca) \to 0$ as $n \to \infty$.

Thus, $d_\HH (S^1, \JRnca)$ and $d_\HH(S^1, \KRnca)$ both tend to zero as $n\to \infty$, hence $\lim_{n\to\infty} \JRnca = \lim_{n\to\infty} \KRnca = S^1$.
\end{proof}

\subsection{Mandelbrots for $c=0$}

The case $c=0$ has been most heavily studied, and for this case, we can obtain more specific results than the general $c \neq 0$.  
Thus, we start by considering the subfamily 
$$\Rna (z) = R_{n, 0, a} (z) =z^n + \frac{a}{z^n},$$ 
and denote the Julia sets by $\KRna$ and $\JRna$.

For the case $c=0$, there is only one free critical orbit.  If $n$ is even, then both critical values have the same image.  If $n$ is odd, then the orbits of these two critical values behave symmetrically under $z \mapsto -z$.  
Thus we may define the Mandelbrot set $\MRna$ as the set of all $a$ such that the free critical orbit does not escape to infinity.  See Figure~\ref{fig:MR0}.

\begin{figure}
\drawfigMRnaexmpnlo
\drawfigMRncaexmpnhi
\caption{\label{fig:MR0}$\MRna$ for $n=6$ and $n=15$}
\end{figure}

We'll keep in mind:

{\bf The Fundamental Trichotomy.} \cite{DevRatPert} 
\textit{
 Let $B_{\infty}$ denote the basin of infinity of a map $\Rna$, with critical values $v^{\pm}$.
\begin{enumerate}
\item if one of $v^{\pm}$ lies in $B_{\infty}$, then $\JRna$ is a Cantor set;
\item or if not, but if one of $\Rna(v^{\pm})$ lies in $B_{\infty}$, then
$\JRna$ is a Cantor set of simple closed curves (i.e., $a$ is in the McMullen domain);
\item or, if the orbit of $v^{\pm}$ is bounded, then $\JRna$ is connected.
\end{enumerate}
}

Hence $\MRna$ is also the connectivity locus for this family.  
We will obtain the following analog of Theorem~\ref{thm:mainMPnc}:

\begin{prop} \label{thm:RM0} Let $S^1_{1/4}$ denote the circle of radius $1/4$ in $\CC$.  Then $$\displaystyle{\lim_{n\to \infty} \MRna = S^1_{1/4}.}$$ \end{prop}

\hide{
Using the lemmas on $K$ and $J$ above, we can establish the geometric limit of each of these classes:


\begin{prop} \label{prop:MR0inS}
For the family $\Rna$, as $n \to \infty$, the Cantor locus tends to $\CC\setminus \DD_{1/4}$,
and the McMullen domain tends to $\DD_{1/4}$.  
Hence  $\lim_{n \to \infty} \MRna \subset S^1_{1/4}$.
\end{prop}

\begin{proof}
If $|a| > 1/4$, then $| v^{\pm} | > 1$. Hence by Lemma~\ref{lem:RcKinD}, for $n$ sufficiently large $v^{\pm} \in B_{\infty}$, thus we have the first statement.  For the second statement, if $|a| < 1/4$, then $|v^{\pm}| < 1$, hence $v^{\pm} \notin B_{\infty}$, but Lemma~\ref{lem:RcDoutK} shows that for $n$ sufficiently large, $\Rna(v^{\pm}) \in B_{\infty}$, thus by the fundamental trichotomy, $a$ is in the McMullen domain for all $n$ sufficiently large.  
\end{proof}
}  

First, we translate the lemmas on $K$ and $J$ above to $\MRna$:

\begin{lem} \label{lem:MRnainA}
For any $\eps \in (0,1/4)$, there is an $N \geq 2$ such that for all $n \geq N$, $\MRna \subset \An{1/4-\eps}{1/4+\eps}$.
\end{lem}

\begin{proof}
Let $\eps \in (0,1/4)$.  Suppose $|a| \geq 1/4 + \eps$. Then 
$$|v^{\pm}| = 2|\sqrt{a}| > \sqrt{1+4\eps} = 1 + \delta$$
for some $\delta > 0$.
By Corollary~\ref{cor:RcKinA}, we know there is an $N \geq 2$ so that for all $n \geq N, \KRna \subset \An{1-\delta/2}{1+\delta/2}.$  Thus, $v^{\pm} \in B_{\infty}$, hence $a \notin \MRna$.  
Hence, $\MRna \subset \DD_{1/4+\eps}$.

Now suppose $|a| \leq 1/4 - \eps$.  Then
$$|v^{\pm}| = 2|\sqrt{a}| < \sqrt{1-4\eps} = 1 - \delta$$
for some $\delta \in (0,1)$.
Then again by Corollary~\ref{cor:RcKinA}, we know there is an $N_1 \geq 2$ so that for all $n \geq N_1, \KRna \subset \An{1-\delta/2}{1+\delta/2}.$  Thus, $v^{\pm}$ are in the trapdoor, i.e., $v^{\pm} \notin B_{\infty}$, but Lemma~\ref{lem:RcDoutK} shows that for all $n > N$ for some $N \geq N_1$,  we have $\Rna(v^{\pm}) \in B_{\infty}$. Thus by the fundamental trichotomy, $a$ is in the McMullen domain for all $n$ sufficiently large, hence not in $\MRna$.  
Thus $\MRna \subset \CC \setminus \cDD_{1/4-\eps}$.
\end{proof}

This lemma implies that  $\inf_{a \in \MRna} d(a, S^1_{1/4}) \to 0$ as $n \to \infty$.  

In order to establish the other distance bound, we will locate specific points in $\MRna$ that will  fill up $S^1_{1/4}$ in the limit.
As in Figure~\ref{fig:MR0}, one can see that for $n \geq 3$, $\MRna$ contains $(n-1)$ baby mandelbrot sets, with one always along the positive real axis, with center tending to $1/4$ as $n \to \infty$, and the others spread evenly along the circle of the same radius.  
We calculate these centers to show existence of a dense subset of $S^1_{1/4}$ in the limit of the $\MRna$'s.  This calculation is a special case of results in \cite{Devrings}. For accessibility, we include here the details needed to support our explicit results without assuming the results in \cite{Devrings}.

For any $n \geq 3$, and each $k \in \{0, 1, \ldots n-2\}$, set 
$$a_n^k = 2^{-2n/(n-1)} e^{i 2\pi k/(n-1)}.$$   Denote the critical values of the map $R_{n,a_n^k}$  by $v^+(a_n^k) = 2 \sqrt{a_n^k}$ and $v^-(a_n^k) = -2 \sqrt{a_n^k}$. 

\begin{lem} \label{lem:RcentersM}
Let $n \geq 3$ and  $k \in \{0, 1, \ldots, n-2\}$.  Then the free critical orbits of $R_{n, a_n^k}$ are finite, hence $a_n^k \in \MRna$.  In particular, if we let $R$ denote $R_{n, a_n^k}$ and let $v^{\pm}$ denote $v^{\pm}(a_n^k)$, 
then we have:
\begin{itemize}
\item if $n$ is even and $k$ is even, then $R(v^-) = v^+$ and $R(v^+) = v^+$;

\item if $n$ is even and $k$ is odd, then $R(v^+) = v^-$ and $R(v^-) = v^-$;

\item if $n$ is odd and $k$ is even, then $R(v^+) = v^+$ and $R(v^-) = v^-$; or

\item if $n$ is odd and $k$ is odd, then  $R(v^+) = v^-$ and $R(v^-) = v^+$.
\end{itemize}
\end{lem}

\begin{proof}
Note 
$$v^+ = 2 \sqrt{a_n^k} = 2 \cdot 2^{-n/(n-1)} e^{i \pi k/(n-1)},$$ and
$$v^- = -2 \sqrt{a_n^k} = 2 \cdot 2^{-n/(n-1)} e^{i \pi k/(n-1)} e^{i \pi}.$$
Define $\sigma \colon \NN \to \{+1, -1\}$ by $\sigma(k) = 1$ if $k$ is even, and $\sigma(k) = -1$ if $k$ is odd.  Then note $e^{i \pi k} = e^{-i \pi k} = \sigma(k)$.  

Now, we calculate $R(v^+)$:
\begin{align*}
R(v^+) &=  (v^+)^n + \frac{(a_n^k)}{(v^+)^n} \\
& = \( 2 \cdot 2^{-n/(n-1)} e^{i \pi k/(n-1)} \)^n \ + \
		\frac{ 2^{-2n/(n-1)}e^{i 2 \pi k/(n-1)}  }{ \( 2 \cdot 2^{-n/(n-1)} e^{i \pi k/(n-1)} \)^n } \\
& = 2^{-n/(n-1)} e^{i \pi kn/(n-1)} \ + \
		 \frac{ 2^{-2n/(n-1)}e^{i 2 \pi k/(n-1)}  }{ 2^{-n/(n-1)} e^{i \pi kn/(n-1)} }  \\
&= 2^{-n/(n-1)} e^{i \pi kn/(n-1)} \ + \ 2^{-n/(n-1)} e^{-i \pi k (n-2)/(n-1)} \\
&=  2^{-n/(n-1)} \( \(e^{i \pi k/(n-1)}\)^n   \ + \  \( e^{-i \pi k/(n-1)} \)^{n-2} \) \\
&= 2^{-n/(n-1)} \( e^{i \pi k/(n-1)} e^{i \pi k} \ + \  \( e^{-i \pi k/(n-1)}\)^{-1} e^{-i \pi k} \) \\
& = 2^{-n/(n-1)}  \( e^{i \pi k/(n-1)}e^{i \pi k} \ + \ e^{i \pi k/(n-1)} e^{-i \pi k}\) \\
&= \sigma(k) \cdot 2 \cdot 2^{-n/(n-1)} e^{i \pi k/(n-1)}. \\
\end{align*}
Thus, if $k$ is even, then $\sigma(k) = 1$, hence $R(v^+) = v^+$. If $k$ is odd, then $\sigma(k) = -1 = e^{i \pi}$, hence $R(v^+) = v^-$.

Now, applying a similar calculation to $v^- = v^+ e^{i \pi}$ yields
\begin{align*}
R(v^-) &= 2^{-n/(n-1)} \( 
 e^{i \pi k/(n-1)}e^{i \pi k} e^{i \pi n} \ + \ e^{i \pi k/(n-1)} e^{-i \pi k} e^{-i \pi n} \)  \\
 & = \sigma(k) \cdot \sigma(n) \cdot 2 \cdot 2^{-n/(n-1)} e^{i \pi k/(n-1)}.
 \end{align*}

Thus $R(v^-) = v^+$ if both $n$ and $k$ are even, or if both $n$ and $k$ are odd, and $R(v^-) = v^-$ if $n$ is even and $k$ is odd, or $n$ is odd and $k$ is even.
\end{proof}

\begin{proof}[Proof of Proposition~\ref{thm:RM0}]
Let $\eps > 0$. First, as in the proof of Theorem~\ref{thm:mainJPnc}, Lemma~\ref{lem:MRnainA} shows that $\inf_{a \in \MRna} d(a, S^1_{1/4}) \to 0$ as $n \to \infty$.

Now, by the above lemma, since the points $a_n^k$ are $n$ points in $\MRna$ spread evenly along a circle of radius $2^{-2n/(n-1)}$, and that radius tends to $1/4$ as $n\to \infty$, we immediately get that for any $s\in S^1_{1/4}, d(s, \MRna) \to 0$ as $n \to \infty$.

Hence, $d_\HH(S^1_{1/4}, \MRna) \to 0$ as $n\to\infty$, thus $\lim_{n\to\infty} \MRna = S^1_{1/4}$.
\end{proof}

This proposition establishes the case $c=0$ of Theorem~\ref{thm:MRnca}.  

\subsection{Mandelbrots for $c\neq0$}

When $c \neq 0 $, the critical orbits of $v^+=c+2\sqrt{a}$ and $v^-=c-2\sqrt{a}$ are both free. This makes the parameter space behavior significantly more complicated.

We denote $\MRnca = \MoneRnca = \{ a \in \CC^* \colon $ at least one free critical orbit of $\Rnca$ is bounded $\}$.  So this is a Mandelbrot living in the slice $\{ c \} \times \CC$ of the parameter space $\CC^* \times \CC^* = \{ (c,a) \colon c, a \in \CC^* \}$.

We let $\MtwoRnca$ denote the $a \in \CC^*$ such that both free critical orbits are bounded.
However, we will focus on $\MoneRnca$, since below we show that requiring both critical orbits to be bounded is too strong of a condition for our purposes:

 \begin{prop}
 If $\abs{c} > 1$, then there is an $N \geq 2$ such that for all $n \geq N$, we have $\MtwoRnca = \emptyset$.
 \end{prop}
 
\begin{proof}
Let $a \in \CC^*$.  Set $w = 2\sqrt{a}$.  

If $\abs{\Arg{w}-\Arg{c}} < \pi$, then since $w \neq 0$, $\abs{c+w} > \abs{c}$.  Hence, $\abs{v^+} = \abs{c+w} = 1+\eps$ for some $\eps > 0$.  Hence there is an $N \geq 2$ such that for all $n \geq N$, the orbit of $v^+$ escapes, thus $a \notin \MtwoRnca$.

On the other hand,  if $\abs{\Arg{w}-\Arg{c}} \geq \pi$, then $\abs{\Arg{-w}-\Arg{c}} < \pi$, hence $\abs{c-w} > \abs{c}$, and again $a \notin \MtwoRnca$.

\end{proof}


Now we turn to establishing that $\lim_{n\to\infty} \MRnca$ is the lima\c con $L_c$.  Theorem~\ref{thm:RcJ} shows us the correct approach:  the lima\c con arises as the set of $c$ for which the critical values have modulus one. 

\begin{prop} \label{prop:Lcformula}
For any $c \in \CC$, $L_c$ is precisely the set of all $a \in \CC$ for which $|c + 2\sqrt{a}| =1$ and/or 
$|c-2\sqrt{a}|=1$.
\end{prop}

\begin{proof}
$\abs{c \pm 2\sqrt{a}} = 1$ iff there is some $\theta$ such that
$c \pm 2\sqrt{a} = - e^{i\theta}$, and we solve for $a$ to get:
\begin{align*}
a & = \left(  \frac{e^{i \theta} + c}{2} \right) ^2 \\
& = \frac{c^2}{4} + \frac{c}{2} e^{i\theta} + \frac{1}{4} e^{2i\theta} \\
& = \left (  \frac{1}{4} + \frac{c}{2}e^{i\theta} + \frac{1}{4}e^{2i\theta}   \right ) + \left(  \frac{c^2-1}{4} \right).
\end{align*}

Now,  suppose $c > 0$.  The polar lima\c con $r=B+A\cos\theta$ has complex equation
$z=\frac{A}{2} + B e^{i\theta} + \frac{A}{2} e^{i2\theta}$, hence $A=1/2$ and $B=c/2$ yields $r=(c+\cos\theta)/2$.

Recall that for $c >0$, we defined $L_c$ as the polar lima\c con $r=(c + \cos\theta)/2$ shifted by $(c^2-1)/4$, i.e.,
$$ L_{c} = \left\{ \left( \frac{1}{4} + \frac{c}{2} e^{i\theta} +  \frac{1}{4}  e^{2i\theta} \right ) 
+ \left( \frac{c^2 -1}{4}  \right) \colon \theta \in [0, 2\pi) \right\}. $$

Now consider $c$ complex, and let $\phi=\Arg{c}$.  Since we defined $L_c := e^{2i\phi} L_{|c|}$, we have
$a \in L_c$ iff $e^{-2i\phi} a \in L_{|c|}$, which means
\begin{align*}
1 & = \abs{|c| \pm 2 \sqrt{a e^{-i2\phi}}} \\
& = \abs{e^{i \phi}}  \abs{|c| \pm 2 \sqrt{a e^{-i2\phi}}} \\
& =  \abs{|c| e^{i \phi} \pm 2 \sqrt{a} } \\
& = \abs{ c \pm 2 \sqrt{a}}. 
\end{align*}

\end{proof}

\begin{figure}
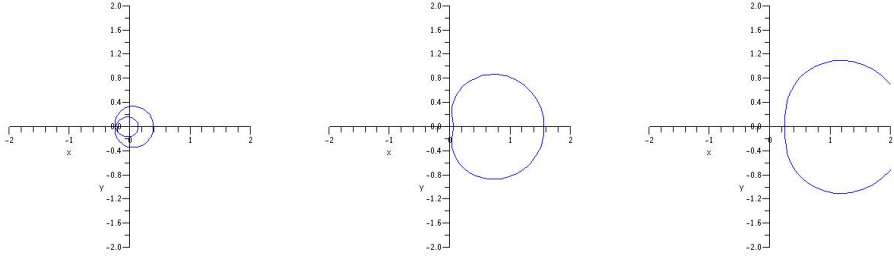

\drawfiglimclo
\drawfiglimcmed
\drawfiglimchi
\caption{\label{fig:Lc}$L_c$ for $c=0.25$; $c=1.5$; and $c=2$ (L to R).}
\end{figure}

\begin{figure}
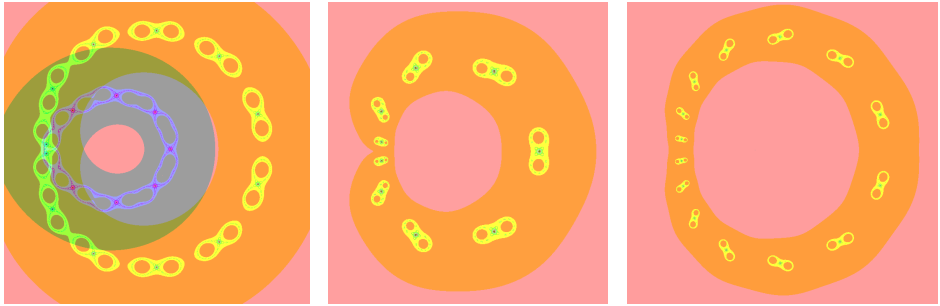

\drawfigMlimclo
\drawfigMlimcmed
\drawfigMlimchi
\caption{\label{fig:MRc}$\MRnca$ for $c=0.25, n=20$; $c=1.5, n=9$; and $c=2, n=14$ (L to R).}
\end{figure}

\begin{figure}
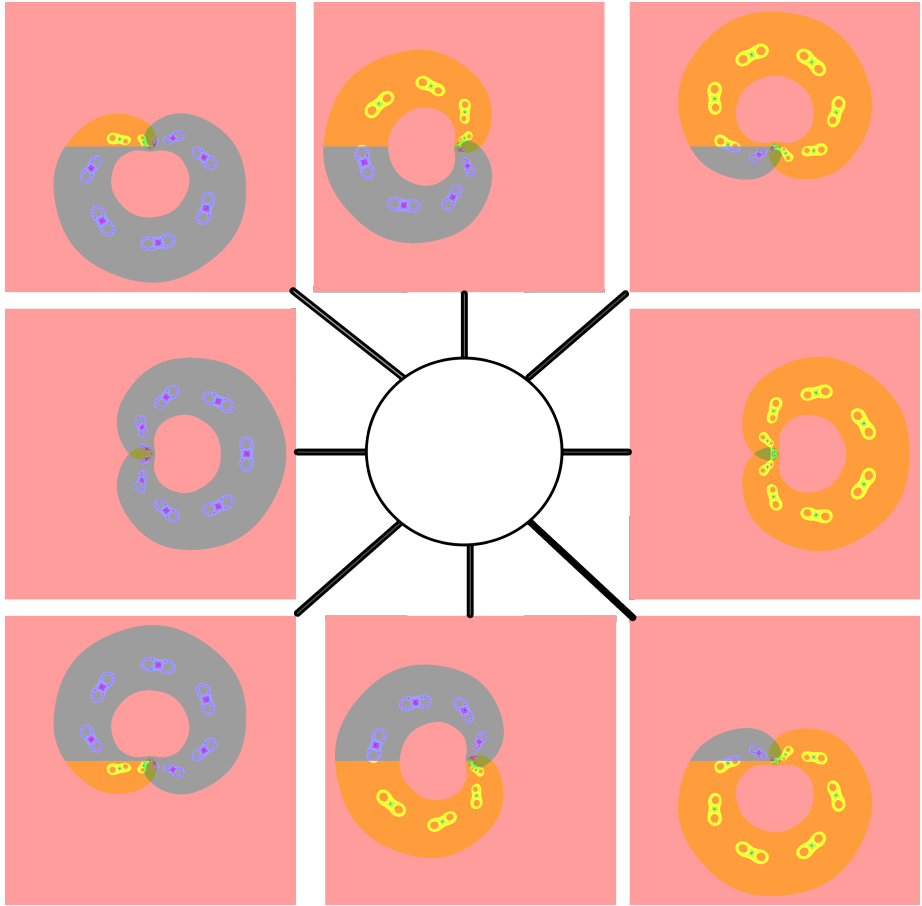

\drawfigcardrot
\caption{\label{fig:cardioidrot} The central circle represents an $S^1$ of $c$-values.  The images about the circle center above are each an $\MRnca$ for $|c|=1$. Each figure is connected by a straight line to its $c$-value. The orange/grey change in coloration is due to a choice of square root.}
\end{figure}

See Figure~\ref{fig:Lc} for images of some lima\c cons $L_c$, and compare to the Mandelbrots $\MRnca$ in Figure~\ref{fig:MRc}.   See Figure~\ref{fig:cardioidrot} for various $\MRnca$ along a circle of common radius $|c|=1$.


First, we will restrict $\MRnca$ to an annulus based on $|c|$.

{\bf Notation.}  Let $c \in \CC$. 
Set
$$l(c)= \frac{(1-\abs{c})^2}{4} \ \text{ and } \ u(c )= \frac{(1+\abs{c})^2}{4},$$
and define
$A_{|c|} =  \overline{\An{l(|c|)}{u(|c|)}}$.

\medskip


\begin{lem} \label{lem:MRncainA}
Let $c \in \CC^*$. Let $\eps > 0$.  Then there is an $N \geq 2$ such that for all $n \geq N$, we have
$$\MRnca \subset \mathcal{N}_\eps(A_{|c|}).$$
\end{lem}

\begin{proof}
First, consider $\CC \setminus \DD_{u(|c|)}$.  Define $\delta > 0$ by
$$
2 \sqrt{ \frac{1}{4} (1+|c|)^2 + \eps } - \abs{c} = 1+ \delta.
$$
Suppose $a$ satisfies $\abs{a} \geq u(c) + \eps$, note $\abs{c} \leq 2\abs{\sqrt{a}}$, and we get
$2\sqrt{|a|} - |c| \geq 1+ \delta$.  Hence by Corollary~\ref{cor:RcKinA}, we see there is an $N_1\geq 2$ (which depends only on $c$ and $\eps$) such that for all $n \geq N_1$, $v^{\pm} \notin \KRnca$, thus $a \notin \MRnca$.  Thus $\MRnca \subset \DD_{u(|c|) + \eps}$ for all $n \geq N_1$.

Next, consider $\DD_{l(|c|)}$.  

First, assume $|c| > 1$.  Define $\delta > 0 $ by
$$
|c| - 2\sqrt{\frac{1}{4} (1-|c|)^2 - \eps } = 1 + \delta.
$$
Then $\abs{a} < l(c)-\eps$ implies $\abs{c} \geq 2\abs{\sqrt{a}} + 1,$
so  $|v^{\pm}|  \geq |c| - 2\sqrt{|a|} \geq 1 + \delta$. Hence by Corollary~\ref{cor:RcKinA}, we see there is an $N\geq N_1$ (which depends only on $c$ and $\eps$) such that for all $n \geq N$, $v^{\pm} \notin \KRnca$, thus $a \notin \MRnca$.  Hence $\MRnca \subset  \mathcal{N}_\eps(A_{|c|})$ 
for all $n \geq N$.

On the other hand, if $\abs{c} < 1$,
define $\delta > 0$ by
$$
|c| + 2{\sqrt{ \frac{1}{4}(1-|c|)^2 - \eps  }} = 1 - \delta.
$$
Then $\abs{a} < l(c) - \eps$ implies
$|v^{\pm}|  \leq \abs{c} + 2\abs{\sqrt{a}} \leq 1 - \delta,$ so by Corollary~\ref{cor:RcKinA}, we see there is an $N\geq N_1$ (which depends only on $c$ and $\eps$) such that for all $n \geq N$, $v^{\pm} \notin \KRnca$, thus $a \notin \MRnca$.  Hence $\MRnca \subset  \mathcal{N}_\eps(A_{|c|})$ for all $n \geq N$.

Finally, note if $\abs{c} =1$, then $l(|c|) = 0$, hence $A_{|c|} = \cDD_{u(|c|)}$ so we are already done.
\end{proof}

Note we get that if $\abs{c} > 1$, then  for large $n$ the complement of this annulus is in the Cantor locus.  
On the other hand, if $\abs{c} < 1$ then for large $n$, $\DD_{l(c)} $ is contained in the McMullen domain and $\DD_{u(c)}$ is in the Cantor locus.  
Finally, if $\abs{c} =1$, we get $\overline{\An{l(c)}{u(c)}}= \cDD$, and outside this disk is the Cantor locus.


\hide{
First,  we restrict the location of $L_c$ to an annulus determined by $\abs{c}$.

\begin{lem} \label{lem:LcinA}
Let $c \in \CC^*$.
%
Then $L_c \subset A_{|c|}$.
\end{lem}

\begin{proof}
We use Proposition~\ref{prop:Lcformula}, and will show that $a \notin A_{|c|}$ implies $a \notin L_c$.

First, suppose that $\abs{a} > u(c) = (1+\abs{c})^2/4.$ Then $\abs{c} \leq 2\abs{\sqrt{a}}$, and
$$\abs{c \pm 2 \sqrt{a}}  \geq  2 \abs{\sqrt{a}} - \abs{c} > 1.$$

Next, suppose $\abs{a} < l(c) = (1-\abs{c})^2/4$.  

If $\abs{c} > 1$, then $\abs{a} < l(c)$ implies $\abs{c} \geq 2\abs{\sqrt{a}} + 1,$ so
$$ \abs{c \pm 2 \sqrt{a}}  \geq \abs{c} - 2 \abs{\sqrt{a}} > 1.$$

On the other hand, if $\abs{c} < 1$, then $\abs{a} < l(c)$ implies
$$\abs{c \pm 2\sqrt{a}} \leq \abs{c} + 2\abs{\sqrt{a}} < 1.$$

Finally, if $\abs{c} =1$ then $\abs{a} > 1$ implies $\abs{2 \sqrt{a}} > 2 = \abs{c} +1$, thus
$$\abs{c \pm 2 \sqrt{a}}  \geq  2 \abs{\sqrt{a}} - \abs{c} > 1.$$
\end{proof}
} 

\begin{lem} \label{lem:MRncainL}
Let $c \in \CC^*$.  Let $\eps > 0$.  Then there is an $N \geq 2$ such that for all $n \geq N$, $\MRnca \subset \mathcal{N}_\eps (L_c)$.
\end{lem}

\begin{proof}
First apply Lemma~\ref{lem:MRncainA}, to choose $N_1 \geq 2$ so that for all $n \geq N_1$, 
$\MRnca \subset \mathcal{N}_\eps(A_{|c|})$.  
Also, note the proof of that lemma combined with Proposition~\ref{prop:Lcformula} implies that $L_c \subset A_{|c|}$.
Hence if $a \notin \mathcal{N}_\eps (L_c)$, there are two cases.  First, if $a \notin \mathcal{N}_\eps(A_{|c|})$, then Lemma~\ref{lem:MRncainA} gives us that for every $n \geq N_1$, $a \notin \MRnca$, and we are done.  Hence, consider 
$${U}_c(\eps) := \overline{\mathcal{N}_\eps(A_{|c|})} \setminus  \mathcal{N}_\eps (L_c).$$
This set is compact, hence the open cover
$\{ B(w, \frac{\eps}{2}) \colon w\in U_{c}(\eps) \}$ has a finite subcover $\{ B_1(\frac{\eps}{2}), \ldots, B_m(\frac{\eps}{2}) \}$.  

Now, if $|c|=1$ then $A_{|c|}=\cDD_{u(|c|)}$, but inspection of the formula for $L_c$ shows that $L_c$ contains the origin (it is the cusp point of the cardioid), hence $U_c(\eps)$ avoids $\DD_{\eps}$ and the covers of balls $B(\frac{\eps}{2})$ avoid $\DD_{\frac{\eps}{2}}$.  Thus, in each ball $B(\frac{\eps}{2})$ the square root function produces two open neighborhoods that are square roots of each $B(\frac{\eps}{2})$.  On the other hand, if $|c| \neq 1$, the annulus $A_{|c|}$ avoids a neighborhood of the origin, hence for $\eps$ sufficiently small, the balls $B(\frac{\eps}{2})$ also avoid a neighborhood of the origin, and the square root function again produces two open neighborhoods as square roots of each $B(\frac{\eps}{2})$.

Hence $\{ c  \pm 2\sqrt{B_k} \}$ is a collection of open neighborhoods, and each $\overline{B_k} \subset \CC \setminus  \mathcal{N}_{\frac{\eps}{2}} (L_c)$, so for each $B_k$ there is a $\delta_k = \delta_k (c, \eps) > 0$ defined by
$$
\delta_k = \min_{w \in \overline{B_k}} \abs{ \abs{c \pm 2\sqrt{w}} -1}.
$$
Hence we may define
$\delta = \min_{1 \leq k \leq m} \delta_k$, and $\delta > 0$.

Now, for this $\delta$, Corollary~\ref{cor:RcKinA} produces an $N \geq N_1$ so that for each $n \geq N$, $\KRnca \subset \mathcal{N}_{\delta}(S^1)$.  But for every $a \in U_c(\eps)$, we have $a \in B_k(\frac{\eps}{2})$ for some $k$, hence $||c \pm 2\sqrt{a}| - 1| > \delta$, thus $v^{\pm} \notin \KRnca$ and so $a \notin \MRnca$.  

Thus, for all $n \geq N$, we have $\MRnca \subset (\CC \setminus U_c(\eps)) \cap  \mathcal{N}_\eps(A_{|c|}) =  \mathcal{N}_\eps(L_{c})$.

\end{proof}

Note this lemma implies that $\inf_{a \in \MRnca} d(a, L_c) \to 0$ as $n\to \infty$. We turn now to the other half of the distance calculation needed.  In the previous section we were able to calculate explicit parameters $a \in \MRna$.  Here we settle for showing the existence of parameters $a \in \MRnca$ which fill up $L_c$ as $n$ grows.

\begin{prop}
For any $c \in \CC^*$, we have $\displaystyle\inf_{a \in L_c} d(a, \MRnca) \to 0$ as $n \to \infty$. 
\end{prop}

\begin{proof}
We seek parameters $a=a_{n,k}(c)$ that solve either equation 
\begin{equation} \label{eqn:critfixed}
c \pm 2\sqrt{a} = a^{1/2n}.
\end{equation}
For in that case, the map has a critical point equal to a critical value, hence has a fixed critical point.  Thus that critical point does not escape, hence  $a \in \MRnca$.  

 
Use $w^{2n} = a$ and transform the above equation to
\begin{equation} \label{eqn:solvecrit}
w^n = \frac{w}{2} - \frac{c}{2}.
\end{equation}

First, some preliminaries.

\begin{figure}
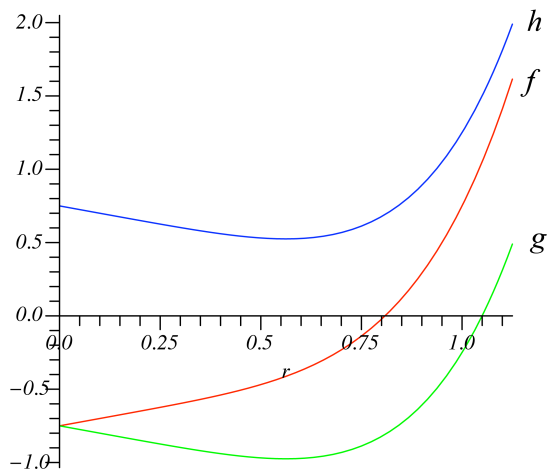

\drawfigfghone
\caption{\label{fig:critsolve1} The functions $f(r) = r^n + \frac{r}{2} - \frac{|c|}{2} $ in red, 
$g(r) = r^n - \frac{r}{2} - \frac{|c|}{2}$ in green, and $h(r) = r^n - \frac{r}{2} + \frac{|c|}{2}$ in blue.  Here we illustrate the case $|c| \geq 1$ with $c=1.5, n=5$.}
\end{figure}
\begin{figure}
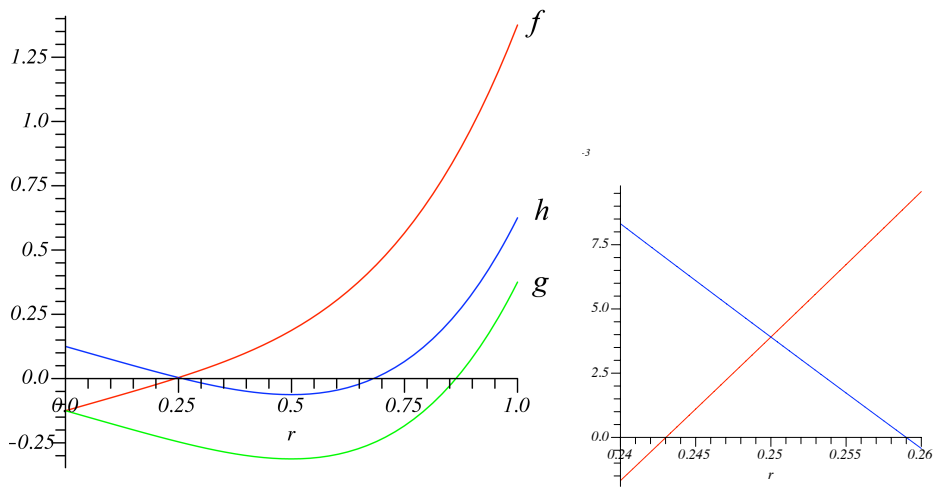

\drawfigfghtwo
\drawfigfghthree
\caption{\label{fig:critsolve2} The functions $f(r) = r^n + \frac{r}{2} - \frac{|c|}{2} $ in red,  
$g(r) = r^n - \frac{r}{2} - \frac{|c|}{2}$ in green, and $h(r) = r^n - \frac{r}{2} + \frac{|c|}{2}$ in blue.  Left: we illustrate the case $|c| < 1$ with $c=0.25, n=4$.  Right: a zoom in near the intersection point of $f$ and $h$, illustrating $u_n < |c| < x_n$.}
\end{figure}

\medskip

\underline{Claim 1:} The function
$$
f(r) = r^n + \frac{r}{2} - \frac{|c|}{2}
$$
has a unique positive root, $r = u_n > 0$.

\underline{Proof:}  $f'(r) = n r^{n-1} + 1/2$ so $f$ has critical points, hence potential optima, at solutions of $r^{n-1} = -1/(2n)$.  Now since $n \in \ZZ^+$ there are two cases.  
(i) If $n$ is even, then $n-1$ is odd, so $f$ has at most one optima (and it is in $\RR^{-}$). Since $f(r) > 0$ for $r \ll 0$ and $r\gg 0$,  $f$ has one minimum.  Since $f(0) = -|c|/2 <0$ and $f(r) >0$ for large $r$, $f$ has one positive root.
(ii)  If $n$ is odd, then $n-1$ is even, so $f$ has no optima, $f$ is monotone increasing.  Since $f(0) < 0$ and $f(r) > 0$ for $r\gg0$, $f$ has a unique positive root.

\medskip

\underline{Claim 2:} The function
$$
g(r) = r^n - \frac{r}{2} - \frac{|c|}{2}.
$$
has a unique positive root, $r = v_n > 0$.

\underline{Proof:} $g'(r) = n r^{n-1} - 1/2$ so $g$ has critical points, hence potential optima, at solutions of $r^{n-1} = 1/(2n)$.  Since $n \in \ZZ^+$ there are two cases.  
(i) If $n$ is even, then $n-1$ is odd, so $g$ has at most one optimum, and it must be in $\RR^+$. Since $g(r) > 0$ for $r \ll 0$ and $r\gg 0$,  $g$ has one minimum.  Since $g(0) = -|c|/2 <0$ and $g'(0)<0$, and $g(r) >0$ for $r\gg 0$, $g$ has one positive root (and one negative root).
(ii)  If $n$ is odd, then $n-1$ is even, so $g$ has up to two optima.  Since $g(r) < 0$ for $r \ll 0$ and $g'(0) <0$, one optimum must be in $\RR^-$.  Since $g(r) >0$ for $r \gg 0$, there must also be one optimum in $\RR^+$.  Since $g(0)<0$, there is a unique root in $\RR^+$.

\medskip

\underline{Claim 3:} $0 < u_n < v_n$.

\underline{Proof:}
Use $\frac{|c|}{2} = (u_n)^n + \frac{u_n}{2} = (v_n)^n - \frac{v_n}{2}$.

\medskip

\underline{Claim 4:} $v_n \to 1$ as $n \to \infty$, and if $|c| \geq 1$, then $u_n \to 1$ as $n \to \infty$. \\(If $|c| < 1$, then $u_n \to |c|$ as $n \to \infty$.)

\underline{Proof:} It is easy to show this by contradiction.

\medskip

\underline{Claim 5:} if $w$ is a solution of Equation~\ref{eqn:solvecrit}, then $u_n \leq |w| \leq v_n$.

\underline{Proof:}
Firstly, $v_n$ is the unique positive 
number satisfying $r^n = \frac{r}{2} + \frac{|c|}{2}$, and for $r > v_n$ we have $r^n >  \frac{r}{2} + \frac{|c|}{2}$.  
But if $w^n =  \frac{w}{2} - \frac{c}{2}$,
then $|w|$ is a positive number satisfying $(|w|)^n = | \frac{w}{2} - \frac{c}{2} | \leq \frac{|w|}{2} + \frac{|c|}{2}$.  So $|w| \leq v_n$.
Secondly, $u_n$ is the unique positive number satisfying $r^n = -\frac{r}{2} + \frac{|c|}{2}$, and if $r < u_n$, then $r^n < -\frac{r}{2} + \frac{|c|}{2}$.  So if $w$ is a solution, then
$|w|^n = \abs{ \frac{w}{2} - \frac{c}{2} } \geq \abs{ \frac{|c|}{2} - \frac{|w|}{2} } \geq \frac{|c|}{2} - \frac{|w|}{2}$, hence $|w| \geq u_n$.

\medskip

The remainder of the proof we divide into two cases.

  \medskip
  
\textbf{Case (1):} assume $|c| \geq 1$.  

\underline{Definition:}
For each $k \in\{ 0, \ldots, n-1\}$, and $\gamma>0$ small, let

$$ S_{n,k} =  
\mathcal{N}_{\gamma}\( \An{u_n}{v_n} \) \cap
\left \{ w \in \CC \colon 
\abs{\Arg{w} - \frac{\Arg{-c}}{n} - \frac{2k\pi}{n} } < \frac{\pi}{n} \right \}.
$$

Each $S_{n,k}$ has lower boundary angle
 $\Theta_{n,k} = \frac{\Arg{-c}}{n} + \frac{(2k-1)\pi}{n}$ and upper boundary angle   $\Theta_{n, k+1}$ 
 and the set of $\overline{S_{n,k}}$ for $k=\{ 0, \ldots, n-1\}$ partitions $\overline{\An{u_n-\gamma}{v_n+\gamma}}$.

\medskip

\underline{Claim 6:} If $|c| \geq 1$, the function
$$
h(r) = r^n - \frac{r}{2} + \frac{|c|}{2}
$$
is strictly positive for $r>0$.  

\textit{Proof:} $h'(r) = nr^{n-1} -\frac{1}{2}$, so  $h$ has critical points, hence potential optima, at solutions of $r^{n-1} = 1/(2n)$.  
(i)  If $n$ is even, then $n-1$ is odd, so $h$ has one critical point. Since $h(r) > 0$ for $r \gg 0$ and $h'(0) <0$, there must be an optimum (minimum) in $\RR^+$, at the positive root of $(\frac{1}{2n})^{\frac{1}{n-1}}$.  We show below that $h$ is positive at that minimum, provided $|c| \geq 1$, thus $h$ is always positive. 
(ii) If $n$ is odd, then $n-1$ is even, so $h$ has two critical points, one positive and one negative.  Again, since  $h(r) > 0$ for $r \gg 0$ and $h'(0) <0$, there must be an optimum (minimum) in $\RR^+$, at the positive root of $(\frac{1}{2n})^{\frac{1}{n-1}}$, and the check below that $h$ is positive at that root guarantees that $h$ is positive for $r > 0$.

Now, assuming $|c| \leq 1$, we have:
\begin{align*}
h \( \(\frac{1}{2n}\)^{\frac{1}{n-1}} \) 
& =  
\( \(\frac{1}{2n}\)^{\frac{1}{n-1}} \)^n - \frac{1}{2} \(\frac{1}{2n}\)^{\frac{1}{n-1}} +  \frac{|c|}{2} \\
 & =  \(\frac{1}{2n}\)^{\frac{n}{n-1}}  - \frac{1}{2} \(\frac{1}{2n}\)^{\frac{1}{n-1}}  + \frac{|c|}{2} \\
 & =  \frac{1}{2n} \(\frac{1}{2n}\)^{\frac{1}{n-1}}  -    \frac{1}{2} \(\frac{1}{2n}\)^{\frac{1}{n-1}}  + \frac{|c|}{2} \\
& =  \(\frac{1}{2n}\)^{\frac{1}{n-1}} \( \frac{1}{2n}  -    \frac{1}{2} \)  + \frac{|c|}{2} \\
& \geq   \frac{1}{2} \( \frac{1}{n} - 1 \) + \frac{|c|}{2} 
 =  \frac{1}{2} \( |c| + \frac{1}{n} - 1 \) \\
& \geq   \frac{1}{2} \( \frac{1}{n} \) > 0.
\end{align*} 

\medskip

\underline{Claim 7:}
For $w$ on the boundary of any $S_{n,k}$,
$$
\frac{|w|}{2} < \abs{w^n + \frac{c}{2}}.
$$

\underline{Proof:}
First consider the  arc boundary segments.  We'll use $\abs{w^n + \frac{c}{2} } \geq \abs{ \frac{|c|}{2} - |w|^n }$. 
On the inner boundary circular arc, we have $|w| < u_n$, hence $|w|^n - \frac{|c|}{2} < -\frac{|w|}{2}$, so $\frac{|c|}{2} - |w|^n > \frac{|w|}{2}$. 
Next, on the outer boundary circular arc, we have $|w| > v_n$, then $|w|^n > \frac{|w|}{2} + \frac{|c|}{2}$, hence $|w|^n - \frac{|c|}{2} > \frac{|w|}{2}$.

Now consider the boundary ray segments.  Suppose $w = re^{i \Theta_k}$ lies on some boundary  ray.  Note $w^n =  r^n e^{i \Arg{-c} } e^{i \pi} = r^n e^{i \Arg{c}}$.  Hence,
$$w^n + \frac{c}{2} 
= \frac{1}{2} \( 2 r^n e^{i \Arg{c}}  + |c| e^{i \Arg{c}} \)
= \frac{1}{2}  e^{i \Arg{c}} \( 2 r^n + |c| \),$$ thus 
$\abs{w^n + \frac{c}{2}} 
= \frac{1}{2} \( 2r^n + |c| \) = r^n + \frac{|c|}{2}.$ 

By claim 6, $r^n + \frac{|c|}{2} > \frac{r}{2} = \frac{|w|}{2}$ is true for any $r>0$, so we have the result.

\medskip

\underline{Claim 8:} Equation~\ref{eqn:solvecrit} has one solution, call it $w_{n,k}(c)$, in each $S_{n,k}$.

\underline{Proof:}
The prior claim yields immediately that, by Rouch\'{e}'s theorem, $w^n - \frac{w}{2} + \frac{c}{2}$ has the same number of zeros as $w^n + \frac{c}{2}$ in each $S_{n,k}$, namely one zero.

\underline{Remainder of proof of Case (1):}  Set $a_{n,k} = (w_{n,k})^{2n}$, then $ c \pm 2\sqrt{a_{n,k}} = w_{n,k}$. 
  Since the $w_{n,k}$ tend to densely fill up $S^1$ as $n\to\infty$, by Proposition~\ref{prop:Lcformula}, we get the $a_{n,k}$ tend to fill up the lima\c con $L_c$ as $n\to\infty$, and since $a_{n,k} \in \MRnca$, we obtain the desired result.
  
  \medskip

\textbf{Case (2):} assume $|c| < 1$.  Refer to Figure~\ref{fig:critsolve2}.
This case is similar in spirit to the above, though a bit simpler.  
Recall $u_n$ is the unique positive root of $f$.  

Note the calculation in the proof of Claim 6 above can be modified to show that for $n$ sufficiently large, $h$ is negative at its positive critical point, hence has two real positive roots $0 < x_n < y_n$, where $x_n$ is close to $|c|$ and $y_n \to 1$ as $n\to \infty$.  

Further, we claim that $u_n < x_n$, and $x_n \to |c|$ as $n \to \infty$.  For, $|c|$ is the positive intersection point of the graphs of $f$ and $h$ (trivial calculation), hence by the shape of the graphs (justified by the calculus done in prior claims), we must have $u_n < |c| < x_n$.  Also, note $f(|c|) = h(|c|) = |c|^n \to 0$ as $n \to \infty$, hence $u_n, x_n \to |c|$ as $n \to \infty$.

Now, for each $k \in\{ 0, \ldots, n-2\}$, and $\gamma>0$ small, let
$$
T_{n,k} =
\mathcal{N}_{\gamma}\( \An{x_n}{y_n} \) \cap
 \left \{ w \in \CC \colon 
\abs{\Arg{w}  - \frac{2k\pi}{n-1} } < \frac{\pi}{n-1} \right \}.
$$
So its boundary angles are $\Phi_k = \pi (2k-1)/{n-1}$ and $\Phi_{k+1}$.

We will show that on the boundary of each $T_{n,k}$, we have
$$
\abs{w^n - \frac{w}{2}} > \frac{|c|}{2}.
$$

First, if $w=re^{i\theta}$ lies on one of the circular arc boundaries, then $r < x_n$ or $r > y_n$, thus $h(r) > 0$, so $r^n -\frac{r}{2} >  \frac{|c|}{2}$.
Thus $\abs{w^n - \frac{w}{2}} \geq \abs{|w|^n - \frac{|w|}{2}} = | r^n  - \frac{r}{2}| > \frac{|c|}{2}$.

Next, if $w=re^{i\phi_k}$ lies on a boundary ray segment,
then $(e^{i \phi_k})^n = e^{i \pi (2k-1) \frac{n}{n-1}} = -e^{i \pi \frac{(2k-1)}{n-1}}$ implies
$$
\abs{w^n - \frac{w}{2}} =  
\abs{-r^n e^{i \pi\frac{(2k-1)}{n-1}} - \frac{r}{2} e^{i \pi\frac{ (2k-1)}{n-1} }}
= \abs{r^n + \frac{r}{2}} = r^n + \frac{r}{2}.$$  
But $x_n > u_n$ implies for all such $r$ we have $f(r) > 0$, hence $ r^n + \frac{r}{2} > \frac{|c|}{2},$ and we are done.

Thus, by Rouch\'{e}'s theorem, the functions $w^n - \frac{w}{2}$ and $w^n - \frac{w}{2} + \frac{c}{2}$ must have the same number of solutions in each $T_{n,k}$.  Since
$
w^n - \frac{w}{2} = w ( w^{n-1} - \frac{1}{2} ) = 0
$ implies $w=0$ or $w=(\frac{1}{2})^{\frac{1}{n-1}}$, we have exactly one solution in each $T_{n,k}$ (note $w=0$ lies in none of the $T's$).

Finally, note that even though $u_n \to |c| < 1$ as $n\to \infty$, we know the solutions $w_{n.k}$ for $n$ large yield $a_{n,k} \in \MRnca$, and we know $\MRnca$ is contained in a neighborhood of $L_c$ which shrinks as $n$ grows, thus the solutions $w_{n,k}$ must satisfy $|w_{n,k}| \to 1$ as $n \to \infty$, and the proof can be finished just as in Case (1).

\end{proof}


Combining this proposition with Lemma~\ref{lem:MRncainA} completes 
the proof of Theorem~\ref{thm:MRnca} in the case $c \neq 0$.


\bibliographystyle{alpha}
\bibliography{multibrot}

 \end{document}